\documentclass[author-year, review]{elsarticle}
\journal{European Journal of Operational Research}

\usepackage[english]{babel}
\usepackage{longtable}
\usepackage{caption}
\usepackage[geometry]{ifsym}
\usepackage{subcaption}
\usepackage[utf8]{inputenc}
\usepackage{multirow}

\usepackage[usenames,dvipsnames]{xcolor}
\usepackage{latexsym,bezier,enumerate,dcolumn}
\usepackage{amsmath,amsthm,amsopn,amstext,amscd,amsfonts,amssymb}
\usepackage{graphics,pstricks,pst-node,lscape,epsfig,fancyhdr,layout}
\usepackage{relsize}
\usepackage{cancel}
\usepackage{marginnote}
\usepackage{setspace}
\usepackage{color}
\usepackage{paracol}
\usepackage{longtable}
\usepackage{arydshln}

\usepackage{tikz}
\usetikzlibrary{automata,arrows,positioning,calc}
\usetikzlibrary{positioning}
\usepackage[colorlinks=true, allcolors=blue]{hyperref}
\usepackage{cleveref}
\usepackage{subcaption}
\usepackage[letterpaper,top=2cm,bottom=2cm,left=3cm,right=3cm,marginparwidth=1.75cm]{geometry}
\usetikzlibrary{patterns}
\usetikzlibrary{shapes.geometric}

\usepackage{blindtext}
\usepackage{anysize}
\usepackage{lmodern}
\usepackage{hyperref}
\usepackage{verbatim}
\hypersetup{%
  colorlinks=true,
  linkcolor={blue!66!black},
  linkbordercolor=red,
  citecolor={blue!66!black}
  }
\definecolor{indigo}{RGB}{75,0,130}
\usepackage[title]{appendix}
\usepackage{arydshln}

\usepackage{algorithm}
\usepackage{algorithmic}
\usepackage{textcomp}

\biboptions{authoryear}

\newcolumntype{d}[1]{D{.}{.}{#1}}

\newcommand{\A}{{\cal A}}
\newcommand{\BB}{{\cal B}}
\newcommand{\C}{{\cal C}}

\newcommand{\E}{{\cal E}}

\newcommand{\G}{{\cal G}}
\renewcommand{\H}{{\cal H}}
\newcommand{\I}{{\cal I}}
\newcommand{\J}{{\cal J}}

\newcommand{\N}{{\cal N}}

\renewcommand{\P}{{\cal P}}
\newcommand{\Q}{{\cal Q}}

\renewcommand{\S}{{\cal S}}
\newcommand{\T}{{\cal T}}

\newcommand{\bea}{\begin{eqnarray}}
\newcommand{\eea}{\end{eqnarray}}

\allowdisplaybreaks[4]

\usepackage{arydshln}
\usepackage[symbol]{footmisc}

 \newcommand{\bcbl}{\color{black}}
 \newcommand{\bcb}{\color{black}}

\begin{document}
\begin{frontmatter}

\title{Multi-horizon optimization for domestic renewable energy system design under uncertainty}

\author[1]{Giovanni Micheli}
\author[2]{Laureano F. Escudero\footnote{Corresponding author.
\newline \hspace*{5.5mm} \textit{Email addresses}: giovanni.micheli@unibg.it (Giovanni Micheli),
laureano.escudero@urjc.es (Laureano F. Escudero),
francesca.maggioni@unibg.it (Francesca Maggioni), bayraksan.1@osu.edu (G{\"u}zin Bayraksan)} }
\author[1]{Francesca Maggioni}
\author[3]{G{\"u}zin Bayraksan}

\address[1]{Department of Management, Information and Production Engineering, University of Bergamo, Italy.}
\address[2]{Area of Statistics and Operations Research, Universidad Rey Juan Carlos, M\'ostoles (Madrid), Spain.}
\address[3]{Department of Integrated Systems Engineering, The Ohio State University, USA.}

\begin{abstract}
In this work we address the challenge of designing optimal domestic renewable energy systems (RES) under multiple sources of uncertainty appearing at different time scales. Long-term uncertainties, such as investment and maintenance costs of different technologies, are combined with short-term uncertainties, including solar radiation, electricity prices, and uncontrolled and static load types.
We formulate the problem as a multistage multi-horizon (MH) stochastic Mixed Integer Linear Programming (MILP) model.
It integrates long-term investment decisions, such as the technology type and capacity of photovoltaic panels and battery energy storage systems, with short-term operational decisions \textcolor{black}{at an hour level}, including energy dispatch, grid exchanges, and supply of a variety of load types.
To ensure robust operation under extreme scenarios, first- and second-order stochastic dominance risk-averse measures are considered preserving the coherency and time consistency of the solution.
Given the computational complexity of solving the stochastic MH MILP for large instances, a rolling horizon-based matheuristic algorithm is considered.
Additionally, various lower-bound strategies are explored, including wait-and-see schemes, expected value approximations, multistage grouping and clustering schemes.
An extensive computational experiment validates the effectiveness of the proposed approach on a case study based on a building complex in South Germany. We tackle models with \textcolor{black}{over 43 million constraints and 12 million binary, 700 hundred integer and 10 million continuous variables, they are solved with up to 0.32\% optimality gap in reasonable computing time, where the value of the stochastic decisions as well as the benefit of the integrated risk-averse measures are quantified.}
\end{abstract}


\end{frontmatter}

\section{Introduction}

Addressing climate change necessitates a coordinated effort to reduce global CO$_2$ emissions, which are predominantly driven by the energy sector.
According to the International Energy Agency (see \cite{WEO22}), the energy industry accounts for approximately 73\% of global greenhouse gas emissions, making it the largest contributor to climate change.
Within this sector, electricity and heat production alone account for roughly 40\% of emissions, totaling 14.7 gigatonnes of CO$_2$ in 2022 (see \cite{GGR23}).
This underscores the urgent need for changes in energy production and consumption practices.

In this context, domestic Renewable Energy Systems (RESs) emerge as a key solution for accelerating the transition to sustainable energy. These systems enable households to generate and use renewable energy, primarily solar, to meet their energy needs, thereby reducing the environmental footprint of energy production.
A typical domestic RES consists of Photovoltaic (PV) panels, Battery Energy Storage Systems (BESS), and digital devices capable of implementing demand response programs by directly controlling household appliances.
Beyond reducing CO$_2$ emissions, domestic RES design contributes to preserving natural resources and ecosystems.
Unlike coal or nuclear power plants,
PV panels operate without water, making them particularly advantageous in regions facing water scarcity.
Moreover, rooftop solar installations minimize land use and habitat disruption, providing an environmentally friendly alternative to large-scale energy projects.

However, the deployment of domestic RES faces several challenges. Sufficient roof space with optimal orientation is essential for effective solar energy capture, which may limit their feasibility for some buildings.
Additionally, the high upfront costs of PV panels and BESS, even when partially mitigated by government incentives, remain a significant barrier for many households.
The intermittent and stochastic nature of solar radiation further complicates the integration of RES into energy systems.
These challenges highlight the need for advanced optimization models to design efficient and cost-effective domestic RES for residential complexes.

There is a broad literature on RES. We refer the readers to \cite{
Zakaria20,
Al-Shahri21,Jasinski23,Ukoba24} for comprehensive reviews on the topic.
Specifically, most of the works presented in the literature are focused on the very short-term, addressing the operations management of a domestic RES that has already been installed; see, for instance, \cite{AlSkaif17,Luo20,OlivieriMcConky20,TaghikhaniZangenehShi22}.
Much fewer works are devoted to the RES design problem.
For example, \cite{Milan12} present a Linear Programming (LP) model for the design of a RES composed of PV panels, a solar thermal collector, and a heat pump over a yearly time horizon, and no uncertainty is considered.
\cite{kucuksarietal_14} formulate a deterministic multistage Mixed Integer Linear Programming (MILP) model  to optimally locate PV panels and design and maintain BESS, and handle the operational uncertainties in the system via a simulation phase that provides updates to the optimization model.
\cite{Mohammadi18} introduce an approach to designing a RES consisting of PV panels, wind turbines, and BESS for a single household over a yearly time horizon; the uncertainty affecting the problem is considered in a post-optimality analysis to validate the solution.

Most studies on RES design adopt deterministic approaches; only a few works in the literature incorporate uncertainty into the modeling framework.
For example, an LP approach is presented in \cite{Zhu21}, where the realizations of wind and solar power output are considered as interval ranges, rather than single numbers, in a fuzzy optimization approach.
\cite{HemmatiSaboori17} propose a stochastic mixed integer non-linear programming model to integrate a BESS into an already installed PV system. \cite{Lauinger16} present a two-stage stochastic LP model to design a RES consisting of PV, fuel cell, heat pump, boiler, solar thermal generators, and BESS. Both works optimize jointly investment decisions on the RES infrastructure and operational decisions on RES components to satisfy an aggregated load under the uncertainty of renewable generation.
The uncertainty of the load is considered in \cite{Aghamohamadi19}, where a robust optimization model is proposed to determine the optimal capacity of a residential PV-battery system while minimizing its operation costs.

The literature on RES also includes some examples of two-stage stochastic MILP models
(see \cite{Schwarz18,CervantesChoobineh18,KeyvandarianSaif24,Beraldi25}).
Specifically, \cite{Schwarz18} develop a two-stage stochastic MILP model to decide, in the first stage, the design of a RES consisting of PV, air-water heat pump, hot water and BESS, and, in the second stage, the system operations in a yearly time horizon with 5 minutes resolution.
Another two-stage stochastic multiperiod MILP model is presented in \cite{CervantesChoobineh18} for PV and BESS design. In this work, an aggregate electricity load is considered and the uncertainty of load, electricity import costs and solar radiation is modeled by means of independent scenarios.
Another example of two-stage stochastic multiperiod MILP model can be found in \cite{KeyvandarianSaif24}, with application to a RES consisting of PV, wind turbines, BESS, and diesel generators.
\bcbl
\cite{Beraldi25} present a two-stage stochastic MILP model for the design of a residential system consisting of PV and BESS, where the available technologies are fixed and the model selects the installed capacities from predefined sets. Multiple types of electricity loads are considered over a multi-period time horizon, and risk aversion is modeled through the Conditional Value-at-Risk (CVaR) measure.
\bcb


As further discussed in \cite{Khezri23} and based on the literature review on RES models, several open challenges can be identified:
\begin{itemize}
\item Integrate multiple time scales within a single framework, combining a medium-to-long-term horizon for RES infrastructure planning with a short-term horizon for system operation and load supply.
\item Provide a more detailed representation of residential electricity demand, disaggregated at the household level and distinguishing between non-controllable loads (e.g., lighting or kitchen appliances) and controllable loads (e.g., air conditioning), \bcbl which can be further classified into elastic and deferrable loads. \bcb
\item Consider accurately both the long-term uncertainty affecting the costs of the RES infrastructure and the short-term uncertainty on weather-based solar radiation and electricity loads.
\item Consider the user discomfort that may be caused by undelivered, \bcbl curtailed or deferred \bcb load in the short-term.
\item Develop solution methods capable of producing high-quality RES configurations with guaranteed performance at a reasonable computational effort.
\end{itemize}

To address these open challenges we deal in this work with the problem of designing a domestic RES under uncertainty by considering
a \textit{multi-horizon stochastic optimization} approach.
This class of problems has been introduced in \cite{Hellemo2012}, see also \cite{Kaut14}, to describe situations that include both
long-term and short-term uncertainty.
The MH approach has been widely utilized in the literature to address real-life problems involving strategic and operational decisions.
To mention just a few:
\cite{WPMHT13} investigate the time consistency property in MH scenario trees;
\cite{ST2015} incorporate short-term uncertainty into long-term energy system models, focusing on power systems with large shares of wind generation;
\cite{SEHT15} introduce a MH stochastic equilibrium model to analyze multi-fuel energy markets, focusing on infrastructure development and renewable energy policies in both perfect and imperfect market structures;
\cite{AA2016} utilize the MH approach to support scheduling decisions for a pumped storage hydro power plant in a liberalized market environment;
\cite{SDPT2016} develop a MH model to plan long-term investments in the European power system considering the uncertainty of the intermittent energy production and the energy demand;
\textcolor{black}{
\cite{Liu17} develop a multistage MH stochastic linear model for generation, storage, and transmission expansion planning, explicitly distinguishing between long-term investment uncertainties and short-term operational variability;
}
\cite{ST2017} investigate through a MH model the influence of policy actions and future energy prices on the cost-optimal development of energy systems in Norway and Sweden;
\cite{cms18} present a MILP model for a MH capacity expansion planning, where a stochastic dominance (SD) functional is considered for risk aversion;
\cite{Maggioni20} extend the definition of the traditional \textit{Expected Value} and \textit{Wait-and-See} problems in stochastic optimization for a MH framework, introducing new measures to quantify the importance of the uncertainty at both strategic and operational levels.
Addressing problems with two time scales is computationally challenging, leading to recent research into efficient solution methods.
\textcolor{black}{\cite{Huang22} introduce a multistage MH stochastic MILP model for power system capacity expansion and propose a nested cross decomposition algorithm to address its computational complexity.}
\cite{Voj23} propose a MH model for generation capacity planning under electricity and natural gas demand uncertainty, leveraging a Benders-type algorithm to achieve near-optimal solutions.
\cite{Zhang24} exploit the block-separable structure of MH stochastic LP models with Benders and Lagrangean decomposition methods. Similarly, \cite{Jacob23} tackle macro-energy systems planning by applying a Benders approach that separates investments from operations and decouples operational time steps via budgeting variables in the master problem.
\cite{Rigaut24} propose a time-block decomposition scheme using recursive Bellman-like equations, demonstrating its effectiveness in a battery management problem, while to address oscillations in Benders schemes, \cite{Zhang24a} introduce a stabilized adaptive method.

In this work, we propose a multistage MH stochastic MILP model to co-optimize long-term decisions in the RES infrastructure to be designed, such as the   capacity of PV panels and BESS, and short-term decisions as well.
The latter ones relate to electricity generation from PV panels, the charge, discharge and energy level of BESS, the electricity exchanges with the power grid, and the supply of different types of load.
Several sources of uncertainty affecting the problem at different time scales are considered: long-term uncertainty comprises the investment and maintenance costs of the different technologies, while short-term uncertainty includes solar radiation, import prices from and export ones to the public power grid, and non-controlled load as well as elastic load types.
To prevent extreme operational scenarios with high discomfort for the users due to undelivered load in the short term,
risk-averse measures are considered.
Given the computational burden of multistage MH stochastic models with risk-averse measures, we also develop a matheuristic algorithm for obtaining feasible solutions and testing their effectiveness by implementing a variety of strategies for obtaining lower bounds.

To sum up, the contributions of our work can be summarized as follows:
{\color{black}
\begin{itemize}\parskip 0.5mm

\item We develop a multistage multi-horizon MILP model for domestic RES design under multiscale uncertainty.
Unlike existing residential RES models, which are predominantly deterministic or two-stage stochastic ones, the proposed formulation integrates long-term investment decisions with detailed short-term operational modeling at the household level, including heterogeneous elastic and deferrable load types within a unified multistage framework.

\item We introduce, for domestic RES planning, coherent and time-consistent risk-averse measures based on first- and second-order stochastic dominance to control extreme short-term discomfort scenarios.
The joint integration of stochastic dominance constraints with a multistage MH residential-scale MILP model constitutes, to the best of our knowledge, a new structural feature in this application domain.

\item We adopt a tactical MH structure in which aggregated operational outcomes, including battery state ones, are explicitly transferred to successor strategic stages.
This inter-stage operational carry-over mechanism tightly couples investment and operational decisions across time, yielding a structurally different formulation from multistage models where operational outcomes do not influence future strategic nodes.

\item The simultaneous presence of cross-scenario risk constraints and inter-stage battery dynamics produces a large-scale strong coupled MILP model that does not satisfy stage-wise or scenario-wise separability.
To address this structural complexity, we adapt and tailor a rolling-horizon-based matheuristic to the specific coupling features of the proposed RES model and complement it with lower-bounding schemes to rigorously assess solution quality.

\item We validate the proposed framework through extensive computational experiments on a real building complex in Southern Germany.
The approach handles instances with over 43 million constraints and 12 million binary variables, achieving optimality gaps as low as 0.32\% within reasonable computational times, and quantifying both the value of multistage stochastic decisions and the impact of risk-averse discomfort control.
\end{itemize}
}

The rest of the paper is organized as follows.
Section \ref{sec:prob} describes the RES design planning problem and provides its formulation as a MH MILP model.
Section \ref{sec:HEUR} describes a constructive matheuristic rolling horizon-based algorithm for providing feasible solutions to the RES design problem.
Section \ref{sec:lb} presents a variety of strategies for obtaining lower bounds.
Section \ref{sec:results} reports the main results of the computational experiments for validating the proposed approach.
Finally, Section \ref{sec:conclu} concludes the paper.
\textcolor{black}{A Supplementary file provides Appendix A with the compact notation of the proposed model, Appendix B with a detailed presentation of established lower bound approaches, and Appendix C describing the procedure adopted for operational scenario generation.}

\vspace*{-0.1in}

\section{Optimal design of a domestic renewable energy system} \label{sec:prob}
\vspace*{-0.05in}
We begin with a \bcbl detailed \bcb description of the domestic RES design problem, followed by its formulation and compact representation. \vspace*{-0.05in}

\subsection{Problem description}
The problem consists of designing a domestic RES by planning investments in PV panels and BESS technologies to satisfy the electricity load of a building complex under several sources of uncertainty.
Two types of decisions must be considered: strategic decisions over a long-term horizon and operational decisions over a short-term horizon.
\bcbl
The strategic decisions are made at the stage set $\E$, where $\E=\{1, \ldots, E\}$, to plan the installation of new capacity for PV and BESS technologies.
Each stage is composed of a number of days $d_e, e\in\E$,  that can differ from one stage to another.
\textcolor{black}{The overall planning horizon is therefore defined implicitly through the stage set and their composition.}
Each day of a stage $e$ is divided into daily periods $t\in\T_e$, where $\T_e$ denotes the corresponding set, represented as blocks of a number, say, $m_t$ of consecutive hours that can vary within the same day. \bcb
The hours of the same daily period are assumed to have the same electricity generation and load.
The management of the operations in the designed RES are made at the daily periods.
In such a MH setting, the two time horizons are affected by several sources of uncertainty: the strategic uncertainty (in the long-term) includes the investment and maintenance costs of the different PV and BESS technologies and their residual values, while the short-term uncertainty comprises solar radiation, electricity prices for selling to and importing from the power grid, load, and operating costs.

In a MH approach to modeling uncertainties, the realization of strategic and operational uncertain parameters is represented by discrete values in scenario trees across different time scales.
\bcbl Specifically, long-term uncertainty is captured through a set of strategic nodes $\N=\{0,1,\ldots,N-1\}$, where $N=|\N|$, organized in a multistage scenario tree, with $\Omega$ denoting the set of strategic scenarios.
On the other hand, operational uncertainty is modeled through a set of operational nodes $\Q_e, e\in\E$ organized within operational subtrees rooted with each strategic node $n$, where $\Pi_{e^n}, e^n\in\E$ denotes the set of scenarios and $e^n$ is the stage to which node $n$ belongs.  \bcb

\begin{figure} [ht!]
	 \begin{center}
     
      \begin{subfigure}[b]{0.49\textwidth}
        \begin{center}
  \resizebox{0.7\textwidth}{0.4\textwidth}{%
	\begin{tikzpicture}[-, >=stealth', auto,
    node distance = 2.5 cm]
 \node[state, red, very thick, fill, scale=0.6] (1) {};
 \node[state, below left=2cm and 3cm of 1, red, very thick, fill,scale=0.6 ] (2)             {} ;
 \node[state, below right=2cm and 3 cm of 1, red, very thick, fill,scale=0.6] (3)              {} ;
 \node[state, below left=2cm and 2cm of 2, red, very thick, fill,scale=0.6 ] (4)              {} ;
\node[state, below right=2cm and 2cm of 2, red, very thick, fill,scale=0.6] (5)              {} ;
 \node[state, below left=2cm and 2cm  of 3, red, very thick, fill,scale=0.6 ] (6)             {}  ;
\node[state, below right=2cm and 2cm of 3, red, very thick, fill,scale=0.6] (7)              {} ;

\node[state, blue,  below=0.5cm of 1, scale = 0.3, diamond, fill] (A12) {A2};
\node[state, blue,  left=0.2cm of A12, scale = 0.3, diamond, fill] (A11) {A1};
\node[state, blue,  right=0.2cm of A12, scale = 0.3, diamond, fill] (A13) {A3};
\node[state, blue,   below=0.3cm of A11, scale = 0.3, diamond, fill] (B11) {B1};
\node[state, blue,   below=0.3cm of A12, scale = 0.3, diamond, fill] (B12) {B2};
\node[state, blue,   below=0.3cm of A13, scale = 0.3, diamond, fill] (B13) {B3};

\node[state, blue,  below=0.5cm of 2, scale = 0.3, diamond, fill] (A22) {A2};
\node[state, blue,  left=0.2cm of A22, scale = 0.3, diamond, fill] (A21) {A1};
\node[state, blue,  right=0.2cm of A22, scale = 0.3, diamond, fill] (A23) {A3};
\node[state, blue,   below=0.3cm of A21, scale = 0.3, diamond, fill] (B21) {B1};
\node[state, blue,   below=0.3cm of A22, scale = 0.3, diamond, fill] (B22) {B2};
\node[state, blue,   below=0.3cm of A23, scale = 0.3, diamond, fill] (B23) {B3};

\node[state, blue,  below=0.5cm of 3, scale = 0.3, diamond, fill] (A32) {A2};
\node[state, blue,  left=0.2cm of A32, scale = 0.3, diamond, fill] (A31) {A1};
\node[state, blue,  right=0.2cm of A32, scale = 0.3, diamond, fill] (A33) {A3};
\node[state, blue,   below=0.3cm of A31, scale = 0.3, diamond, fill] (B31) {B1};
\node[state, blue,   below=0.3cm of A32, scale = 0.3, diamond, fill] (B32) {B2};
\node[state, blue,   below=0.3cm of A33, scale = 0.3, diamond, fill] (B33) {B3};

\node[state, blue,  below=0.5cm of 4, scale = 0.3, diamond, fill] (A42) {A2};
\node[state, blue,  left=0.2cm of A42, scale = 0.3, diamond, fill] (A41) {A1};
\node[state, blue,  right=0.2cm of A42, scale = 0.3, diamond, fill] (A43) {A3};
\node[state, blue,   below=0.3cm of A41, scale = 0.3, diamond, fill] (B41) {B1};
\node[state, blue,   below=0.3cm of A42, scale = 0.3, diamond, fill] (B42) {B2};
\node[state, blue,   below=0.3cm of A43, scale = 0.3, diamond, fill] (B43) {B3};

\node[state, blue,  below=0.5cm of 5, scale = 0.3, diamond, fill] (A52) {A2};
\node[state, blue,  left=0.2cm of A52, scale = 0.3, diamond, fill] (A51) {A1};
\node[state, blue,  right=0.2cm of A52, scale = 0.3, diamond, fill] (A53) {A3};
\node[state, blue,   below=0.3cm of A51, scale = 0.3, diamond, fill] (B51) {B1};
\node[state, blue,   below=0.3cm of A52, scale = 0.3, diamond, fill] (B52) {B2};
\node[state, blue,   below=0.3cm of A53, scale = 0.3, diamond, fill] (B53) {B3};

\node[state, blue,  below=0.5cm of 6, scale = 0.3, diamond, fill] (A62) {A2};
\node[state, blue,  left=0.2cm of A62, scale = 0.3, diamond, fill] (A61) {A1};
\node[state, blue,  right=0.2cm of A62, scale = 0.3, diamond, fill] (A63) {A3};
\node[state, blue,   below=0.3cm of A61, scale = 0.3, diamond, fill] (B61) {B1};
\node[state, blue,   below=0.3cm of A62, scale = 0.3, diamond, fill] (B62) {B2};
\node[state, blue,   below=0.3cm of A63, scale = 0.3, diamond, fill] (B63) {B3};

\node[state, blue,  below=0.5cm of 7, scale = 0.3, diamond, fill] (A72) {A2};
\node[state, blue,  left=0.2cm of A72, scale = 0.3, diamond, fill] (A71) {A1};
\node[state, blue,  right=0.2cm of A72, scale = 0.3, diamond, fill] (A73) {A3};
\node[state, blue,   below=0.3cm of A71, scale = 0.3, diamond, fill] (B71) {B1};
\node[state, blue,   below=0.3cm of A72, scale = 0.3, diamond, fill] (B72) {B2};
\node[state, blue,   below=0.3cm of A73, scale = 0.3, diamond, fill] (B73) {B3};

	\path
	(1) edge[red, very thick, fill] (2)
    (1) edge[red, very thick, fill] (3)
    	(2) edge[red, very thick, fill] (4)
    (2) edge[red, very thick, fill] (5)
    	(3) edge[red, very thick, fill] (6)
    (3) edge[red, very thick, fill] (7)
     (1) edge[blue] (A11)
    (A11) edge[blue] (B11)
    (1) edge[blue]   (A12)
    (A12) edge[blue] (B12)
    (1) edge[blue]   (A13)
    (A13) edge[blue] (B13)
  
    (2) edge[blue] (A21)
    (A21) edge[blue] (B21)
    (2) edge[blue]   (A22)
    (A22) edge[blue] (B22)
    (2) edge[blue]   (A23)
    (A23) edge[blue] (B23)

    (3) edge[blue] (A31)
    (A31) edge[blue] (B31)
    (3) edge[blue]   (A32)
    (A32) edge[blue] (B32)
    (3) edge[blue]   (A33)
    (A33) edge[blue] (B33)

    (4) edge[blue] (A41)
    (A41) edge[blue] (B41)
    (4) edge[blue]   (A42)
    (A42) edge[blue] (B42)
    (4) edge[blue]   (A43)
    (A43) edge[blue] (B43)

    (5) edge[blue] (A51)
    (A51) edge[blue] (B51)
    (5) edge[blue]   (A52)
    (A52) edge[blue] (B52)
    (5) edge[blue]   (A53)
    (A53) edge[blue] (B53)

    (6) edge[blue] (A61)
    (A61) edge[blue] (B61)
    (6) edge[blue]   (A62)
    (A62) edge[blue] (B62)
    (6) edge[blue]   (A63)
    (A63) edge[blue] (B63) 

    (7) edge[blue] (A71)
    (A71) edge[blue] (B71)
    (7) edge[blue]   (A72)
    (A72) edge[blue] (B72)
    (7) edge[blue]   (A73)
    (A73) edge[blue] (B73) ;
	\end{tikzpicture}\bigskip \bigskip 
    }
  \end{center}
  \caption {(a) A \textit{traditional} multi-horizon scenario tree. Terminal leaves of operational sub-trees are not linked to the immediate successor investment and operational nodes.}
  \label{MH_Tree_a}
  \end{subfigure}
  \begin{subfigure}[b]{0.49\textwidth}
	 \begin{center}
  \resizebox{0.7\textwidth}{0.5\textwidth}{%
	\begin{tikzpicture}[-, >=stealth', auto,
    node distance = 2.5 cm]
	\node[state, red, very thick, fill, scale=0.6] (1) {};
 \node[state, below left=3.5cm and 3cm of 1, red, very thick, fill,scale=0.6 ] (2) {}             ;
 \node[state, below right=3.5cm and 3 cm of 1, red, very thick, fill,scale=0.6] (3) {}              ;
 \node[state, below left=3.5cm and 2cm of 2, red, very thick, fill,scale=0.6 ] (4) {}              ;
\node[state, below right=3.5cm and 2cm of 2, red, very thick, fill,scale=0.6] (5) {}               ;
 \node[state, below left=3.5cm and 2cm  of 3, red, very thick, fill,scale=0.6 ] (6) {}               ;
\node[state, below right=3.5cm and 2cm of 3, red, very thick, fill,scale=0.6] (7) {}               ;

\node[state, blue,  below=0.5cm of 1, scale = 0.3, diamond, fill] (A12) {};
\node[state, blue,  left=0.2cm of A12, scale = 0.3, diamond, fill] (A11) {};
\node[state, blue,  right=0.2cm of A12, scale = 0.3, diamond, fill] (A13) {};
\node[state, blue,   below=0.3cm of A11, scale = 0.3, diamond, fill] (B11) {};
\node[state, blue,   below=0.3cm of A12, scale = 0.3, diamond, fill] (B12) {};
\node[state, blue,   below=0.3cm of A13, scale = 0.3, diamond, fill] (B13) {};

\node[state, blue,  below=0.5cm of 2, scale = 0.3, diamond, fill] (A22) {};
\node[state, blue,  left=0.2cm of A22, scale = 0.3, diamond, fill] (A21) {};
\node[state, blue,  right=0.2cm of A22, scale = 0.3, diamond, fill] (A23) {};
\node[state, blue,   below=0.3cm of A21, scale = 0.3, diamond, fill] (B21) {};
\node[state, blue,   below=0.3cm of A22, scale = 0.3, diamond, fill] (B22) {};
\node[state, blue,   below=0.3cm of A23, scale = 0.3, diamond, fill] (B23) {};

\node[state, blue,  below=0.5cm of 3, scale = 0.3, diamond, fill] (A32) {};
\node[state, blue,  left=0.2cm of A32, scale = 0.3, diamond, fill] (A31) {};
\node[state, blue,  right=0.2cm of A32, scale = 0.3, diamond, fill] (A33) {};
\node[state, blue,   below=0.3cm of A31, scale = 0.3, diamond, fill] (B31) {};
\node[state, blue,   below=0.3cm of A32, scale = 0.3, diamond, fill] (B32) {};
\node[state, blue,   below=0.3cm of A33, scale = 0.3, diamond, fill] (B33) {};

\node[state, blue,  below=0.5cm of 4, scale = 0.3, diamond, fill] (A42) {};
\node[state, blue,  left=0.2cm of A42, scale = 0.3, diamond, fill] (A41) {};
\node[state, blue,  right=0.2cm of A42, scale = 0.3, diamond, fill] (A43) {};
\node[state, blue,   below=0.3cm of A41, scale = 0.3, diamond, fill] (B41) {};
\node[state, blue,   below=0.3cm of A42, scale = 0.3, diamond, fill] (B42) {};
\node[state, blue,   below=0.3cm of A43, scale = 0.3, diamond, fill] (B43) {};

\node[state, blue,  below=0.5cm of 5, scale = 0.3, diamond, fill] (A52) {};
\node[state, blue,  left=0.2cm of A52, scale = 0.3, diamond, fill] (A51) {};
\node[state, blue,  right=0.2cm of A52, scale = 0.3, diamond, fill] (A53) {};
\node[state, blue,   below=0.3cm of A51, scale = 0.3, diamond, fill] (B51) {};
\node[state, blue,   below=0.3cm of A52, scale = 0.3, diamond, fill] (B52) {};
\node[state, blue,   below=0.3cm of A53, scale = 0.3, diamond, fill] (B53) {};

\node[state, blue,  below=0.5cm of 6, scale = 0.3, diamond, fill] (A62) {};
\node[state, blue,  left=0.2cm of A62, scale = 0.3, diamond, fill] (A61) {};
\node[state, blue,  right=0.2cm of A62, scale = 0.3, diamond, fill] (A63) {};
\node[state, blue,   below=0.3cm of A61, scale = 0.3, diamond, fill] (B61) {};
\node[state, blue,   below=0.3cm of A62, scale = 0.3, diamond, fill] (B62) {};
\node[state, blue,   below=0.3cm of A63, scale = 0.3, diamond, fill] (B63) {};

\node[state, blue,  below=0.5cm of 7, scale = 0.3, diamond, fill] (A72) {};
\node[state, blue,  left=0.2cm of A72, scale = 0.3, diamond, fill] (A71) {};
\node[state, blue,  right=0.2cm of A72, scale = 0.3, diamond, fill] (A73) {};
\node[state, blue,   below=0.3cm of A71, scale = 0.3, diamond, fill] (B71) {};
\node[state, blue,   below=0.3cm of A72, scale = 0.3, diamond, fill] (B72) {};
\node[state, blue,   below=0.3cm of A73, scale = 0.3, diamond, fill] (B73) {};

 \node[state, below=0.3cm of B12, red, very thick, fill, scale=0.4,rectangle] (O1) {};
 \node[state, below=0.3cm of B22, red, very thick, fill,scale=0.4,rectangle ] (O2) {}              ;
 \node[state, below=0.3cm of B32, red, very thick, fill,scale=0.4,rectangle] (O3) {}               ;

	\path
     (1) edge[blue] (A11)
    (A11) edge[blue] (B11)
    (1) edge[blue]   (A12)
    (A12) edge[blue] (B12)
    (1) edge[blue]   (A13)
    (A13) edge[blue] (B13)
  
    (2) edge[blue] (A21)
    (A21) edge[blue] (B21)
    (2) edge[blue]   (A22)
    (A22) edge[blue] (B22)
    (2) edge[blue]   (A23)
    (A23) edge[blue] (B23)

    (3) edge[blue] (A31)
    (A31) edge[blue] (B31)
    (3) edge[blue]   (A32)
    (A32) edge[blue] (B32)
    (3) edge[blue]   (A33)
    (A33) edge[blue] (B33)

    (4) edge[blue] (A41)
    (A41) edge[blue] (B41)
    (4) edge[blue]   (A42)
    (A42) edge[blue] (B42)
    (4) edge[blue]   (A43)
    (A43) edge[blue] (B43)

    (5) edge[blue] (A51)
    (A51) edge[blue] (B51)
    (5) edge[blue]   (A52)
    (A52) edge[blue] (B52)
    (5) edge[blue]   (A53)
    (A53) edge[blue] (B53)

    (6) edge[blue] (A61)
    (A61) edge[blue] (B61)
    (6) edge[blue]   (A62)
    (A62) edge[blue] (B62)
    (6) edge[blue]   (A63)
    (A63) edge[blue] (B63) 

    (7) edge[blue] (A71)
    (A71) edge[blue] (B71)
    (7) edge[blue]   (A72)
    (A72) edge[blue] (B72)
    (7) edge[blue]   (A73)
    (A73) edge[blue] (B73)

    (B11) edge[blue] (O1)
    (B12) edge[blue] (O1)
    (B13) edge[blue] (O1)
    
    (B21) edge[blue] (O2)
    (B22) edge[blue] (O2)
    (B23) edge[blue] (O2) 
    
    (B31) edge[blue] (O3)
    (B32) edge[blue] (O3)
    (B33) edge[blue] (O3)

    (O1) edge[red, very thick, fill] (2)
    (O1) edge[red, very thick, fill] (3)
    (O2) edge[red, very thick, fill] (4)
    (O2) edge[red, very thick, fill] (5)
    (O3) edge[red, very thick, fill] (6)
    (O3) edge[red, very thick, fill] (7) ;
	\end{tikzpicture}
    }
  \end{center}
  \caption {(b) A \textit{tactical} multi-horizon scenario tree. Terminal leaves of each operational sub-tree are connected to the immediate successor investment node, representing the transfer of aggregated information on short-term operations to the immediate successor investment node.}
  \label{MH_Tree_b}
  \end{subfigure}
  \caption{Comparison between a \textit{traditional} and a \textit{tactical} multi-horizon scenario tree with 3 strategic stages, 7 strategic nodes, and, for each strategic node, operational sub-trees consisting of 2 operational periods and 3 operational scenarios each.}
  \label{MH_Tree}
  \end{center}
  \end{figure}

Figure \ref{MH_Tree_a} illustrates an example of a MH scenario tree with three stages and seven strategic nodes, represented by red circles. Each stage corresponds to a time interval within the long-term horizon where new information becomes available and decisions are made. For each stage $e\in\E$, there is a corresponding set of strategic nodes $\N_e\subset\N$, such that each node $n\in\N_e$ represents a specific decision or event.
\bcbl The strategic node at the first stage, $n=0$, is known as the root, while those at the terminal stage, $n \in \N_E$, are called leaves.
Every strategic node, except the root one, is linked to a unique node in the preceding stage referred to as its ancestor $a(n)$. A strategic scenario $\omega\in\Omega$ is defined as a path through the sequence of strategic nodes from the root to a leaf, representing the realization of strategic uncertain parameters across the long-term horizon.
\bcb
At the operational level, short-term uncertainty is modeled independently of the strategic one. In the MH scenario tree shown in Fig. \ref{MH_Tree_a}, short-term uncertainty is represented by two-stage multi-period scenario subtrees rooted with each strategic node.
Each operational subtree contains six nodes, corresponding to three scenarios and two daily periods.
It is important to notice that the terminal leaves of these operational subtrees are not connected to subsequent strategic or operational nodes.
This representation of uncertainty is referred to as the traditional MH scenario tree, as it was introduced in the seminal works of \cite{Hellemo2012} and \cite{Kaut14}.
On the other hand, in this work we consider the \textit{tactical} representation introduced in \cite{forestry20},
\bcbl where the terminal leaves of each operational subtree, rooted at a given strategic node, are connected to the corresponding immediate successor investment node, as shown in Fig. \ref{MH_Tree_b}.
\bcb
Such a representation allows the transfer of aggregated information on short-term operations to the successor investment and operational nodes, thus providing a better representation of long-term dynamics of BESS.

In line with \cite{Milan12,Haq19}, \bcbl PV technology is modeled by assuming that each model type $i\in\I$ consists of identical panels with nominal capacity $F_i$.
The maximum hourly electricity output of a single PV panel in operational node $q\in\Q^n, \, n\in\N$ is given by $(\tilde{z}^R)_{i,q} F_i$, where $(\tilde{z}^R)_{i,q}$ represents the fraction of nominal capacity available due to solar radiation.
At the strategic level, each PV technology is associated with a fixed preparation cost $C1_i^n$,
a unit installation cost $C2_i^n$, and a unit maintenance cost $C3_i^n$, all modeled as strategic stochastic parameters.
Additionally, each PV technology has an uncertain unit residual value $V2_i^n, \, n\in\N_E$ at the end of the planning horizon.
At the operational level, the electricity demand of the domestic building complex can be supplied by PV generation during day-time periods, depending on the installed capacity and the stochastic availability of solar radiation. Generation is allowed only in a subset of daily periods $\mathcal{T}_e^R \subseteq \mathcal{T}e$ at stage $e\in\E$ and incurs stochastic operating costs $C^R_{i,q}$ for each operational node $q\in\Q_e$.
\bcb

To mitigate the variability of the solar power production,  investments in PV panels can be coupled with investments in different BESS technologies $b \in \BB$.
Similar to PV modeling, each BESS technology consists of identical battery units.
\bcbl At the strategic level, each battery technology $b \in \BB$ is characterized by a fixed preparation costs $C4_b^n$, a unit installation cost $C5_b^n$, a unit maintenance cost $C6_b^n$, and a residual value $V5_b^n$ at terminal nodes $n\in\N_E$.
At the operational level, BESS are defined by a nominal capacity $(k')_b$, an energy loss factor $(f^l)_{b,e}$ in stage $e$, charge/discharge depths fractions $\rho_{b,e}^+, \, \rho^-_{b.e}$ and an operational cost $\epsilon_b$ per battery unit. \bcb
If RES production and storage are insufficient to meet the electricity demand, energy can be imported from the power grid
\bcbl in each operational node $q\in\Q_e, \, e\in\E$, incurring electricity import costs $C_q^G$, which are typically much higher than the costs of RES production.
On the other hand, excess generation can be exported to the power grid at price $\P_q$. \bcb

\bcbl
The building complex consists of households and commercial offices, with different electricity load profiles that can be classified into two categories: {\it non-controlled} loads and {\it controlled} loads.
The first group includes all types of consumption that cannot be controlled by the RES central operator and are denoted by $L_q$, for $q\in\Q_e, \, e\in\E$.
These loads must necessarily be satisfied.
Examples include lighting, kitchen appliances, online computing, and communication devices.
Notice that this type of loads are treated as exogenous demand profiles that must always be
satisfied in each operational scenario.

On the other hand, controlled loads can be managed by a central RES or home management system.
This type of loads is explicitly modeled through decision variables that allow shifting
or curtailing part of the demand within predefined flexibility limits.
The integration of non-controllable and controllable types is performed in the operational part of the RES model to ensure feasibility with respect to the non-controllable loads while exploiting the flexibility of controllable loads to reduce costs.
The controlled loads are represented by the set $\J$ and are characterized by mechanisms that allow modulation of their operation in terms of intensity (how much power they consume) and time scheduling (when they operate).
This distinction divides them into two categories:  {\it elastic} loads, indexed by $\J_1 \subseteq \J$, and {\it deferrable} loads, indexed by $\J_2 \subseteq \J$.

Elastic loads can be partially curtailed, since they are loads whose energy consumption can change, and cause user discomfort.
A typical example includes Heating, Ventilation, and Air Conditioning (HVAC) systems (e.g., heating to a lower temperature).
For each elastic load type $j\in\J_1$, a reference consumption level $(L^1)_{j,q}$ is defined for operational node $q\in\Q_e$.
This setpoint represents the target consumption over a subset of daily periods $(\mathcal{T}^1)_{j,e} \subseteq \mathcal{T}_e$ at stage $e$.
Negative deviations from the setpoint are allowed up to a maximum curtailment value $(\overline{\Delta}\ell^1)_{j,t}$ for each daily period $t \in \mathcal{T}_e$, inducing user discomfort as described below.

Deferrable controlled loads, on the other hand, are loads whose operation can be shifted in time (within a predefined time window),
but whose total energy requirement is fixed. For instance, dishwashers can be run at night.
For each deferrable load type $j\in\J_2$, a reference starting period $\tau_{j,e}$ is defined, together with an average hourly electricity requirement $(L^2)_{j,e}$ to be supplied over a consecutive number of hours $(m^1)_{j,e}$ at stage $e$.
The admissible supply periods belong to a time window $(\mathcal{T}^2)_{j,e} \subseteq \mathcal{T}_e$.
Anticipating or delaying the load supply relative to the reference period $\tau_{j,e}$ induces user discomfort.
Deferrable loads include, for example, electricity load for using the washing machine, dryer, dishwasher, oven, iron, battery of electricity vehicles, and vacuum cleaner.
\bcb
Some pairs of deferrable load types cannot be satisfied simultaneously in the same daily period.
These incompatibilities are represented by the set $\mathcal{H}_1 \subseteq \mathcal{J}_2 \times \mathcal{J}_2$, where $(j,j^\prime)\in\mathcal{H}_1$ denotes mutually incompatible load types.
Other pairs are subject to latency-based precedence relationships.
These are represented by the set $\mathcal{H}_2 \subseteq \mathcal{J}_2 \times \mathcal{J}_2$, where $(j,j')\in\mathcal{H}_2$ indicates that load type $j$ must precede load type $j^\prime$ by at least $\overline{\Delta}t_{j,j^\prime}$ daily periods.

The user discomfort associated with elastic load type $j\in\J_1$ is modeled through a penalization parameter $(D^1)_{j,t_q}$, which represents the discomfort incurred per unit of curtailed load in daily period $t_q \in (\T^1)_{j,e}$ at operational node $q\in\Q_e$.
For deferrable load type $j\in\J_2$, discomfort arises from shifting the supply away from the reference starting period $\tau_{j,e}$.
This is quantified by a penalization parameter $(D^2)_{j,t_q}$ for daily period $t_q \in (\T^2)_{j,e}$ at operational node $q\in\Q_e$, $e\in\E$.
In the risk-neutral formulation, the total discomfort in each operational scenario $\pi \in \Pi_e$ is upper bounded by a threshold parameter $\hat{D}_e$.
For each scenario $\pi$, the set $\A_\pi \subseteq \Q_e$ denotes the operational nodes belonging to that scenario.
In the risk-averse formulation, a set of policy profiles $\P_e$ is considered at each strategic stage.
Each profile $p \in \P_e$ specifies a discomfort threshold $\overline{D}^p$.
First-order stochastic dominance constraints bound the probability of exceeding this threshold through a parameter $\overline{\eta}^p$.
Second-order stochastic dominance constraints limit both the maximum excess fraction $\overline{s}^p$ of $\overline{D}^p$ allowed in each scenario and the expected excess fraction $\overline{\overline{s}}^p$ across scenarios.

In this framework, the decision maker aims to design domestic RES infrastructures capable of meeting the electricity demand of the building complex at minimum total investment and operational cost.
At the strategic level, decisions concern investments in PV and BESS technologies over the long-term planning horizon, subject to technical and budget constraints.
\bcbl
For each strategic node $n\in\N$, the binary variable $x_i^n$ indicates whether PV technology $i$ has been installed \textit{by} node $n$, while $\tilde{x}_i^n \in \mathbb{R}^+$ denotes the total number of PV panels of type $i$ installed \textit{by} node $n$.
The binary variable $\alpha_i^n$ takes value 1 if new PV panels of technology $i$ are installed \textit{at} node $n$.
Similarly, $(x^\prime)_b^n$ indicates whether battery technology $b$ has been installed \textit{by} node $n$, and
$(\tilde{x}^\prime)_b^n \in \mathbb{N}$ represents the total number of battery units of type $b$ installed \textit{by} node $n$.
The binary variable $\beta_b^n$ takes value if new battery units of technology $b$ are installed \textit{at} node $n$.

At the operational level, the problem consists of optimizing electricity generation, load supply, storage operation, and
power exchanges with the public power grid over the short-term horizon.
For each operational node $q \in \Q_{e^n}$ and strategic node $n\in\N$, the variable $(z^R)^n_{i,q} \in \mathbb{R}^+$ denotes the hourly electricity generated from PV panels of type $i$, and $(z^G)^n_{q} \in \mathbb{R}^+$ denotes the hourly electricity imported from the grid.
The storage dynamics are captured by $y^n_{b,q} \in \mathbb{R}^+$, representing the energy stored in battery type $b$ at the end of daily period $t_q$, together with $(y^+)^n_{b,q} \in \mathbb{R}^+$ and $(y^-)^n_{b,q} \in \mathbb{R}^+$, representing the hourly charge and discharge of battery type $b$, respectively.
For elastic loads, $\Delta \ell^n_{j,q} \in \mathbb{R}^+$ represents the deviation from the setpoint $(L^1)_{j,q}$ of load type $j$ in node $q$, while for deferrable loads, the binary variable $\delta_{j,q}^n$ equals 1 if load type $j$ starts at daily period $t_q$.
Finally, $s^{p,\pi} \in \mathbb{R}^+$ denotes the discomfort excess over the threshold $\overline{s}^p \overline{D}^p$ for profile $p$ in operational scenario $\pi$, and the binary variable $\eta^{p,\pi}$ equals 1 if the discomfort threshold is exceeded.
\bcb

\subsection{Problem formulation}\label{sec:model}

Optimal decisions about the domestic RES design and its short-term management are determined by minimizing the objective function \eqref{of}, defined as the sum of the expected strategic and operational costs of the RES design, subject to several types of constraints:
\begin{itemize}
    \item Strategic constraints \eqref{stra-cons} related to the RES design and its budget.
    \item Operational constraints \eqref{stoch-2nd-cons} modeling the short-term operations of PV and BESS technologies.
    \item Operational constraints \eqref{stoch-2nd-delta} related to the supply of the loads.
    \item Tactical constraints \eqref{tac-battery_flow_balance} modeling the BESS flow balance between consecutive stages.
    \item The maximum user discomfort constraint \eqref{rn}, which bounds the expected discomfort caused to the users by load curtailment and shifting with respect to setpoints and reference daily periods.
    \item The risk-averse discomfort constraint system \eqref{sd}, which protects the solution against extreme operational scenarios in terms of user discomfort.
\end{itemize}

We now formulate the objective function and the constraint system.
\bcbl The details of the notation used in this section are provided in Tables A.1--A.9 in Appendix A of the Supplementary file. \bcb

\bigskip

\noindent
\begin{subequations}\label{of}
\textbf{Objective function} \eqref{of}
\begin{align}
 z^* \, =  \,
 \min  \  & \sum_{n\in\N}w^n \sum_{i\in\I}
           \bigl ( C1_i^n ( x_i^n - x_i^{a(n)} ) + C2_i^n ( \tilde{x}_i^n - \tilde{x}_i^{a(n)} ) + C3_i^n \tilde{x}_i^n  \bigr ) \label{of-0}\\
+ & \sum_{n\in\N} w^n \sum_{b\in\BB}
            \bigl ( C4_b^n ( (x')_b^n - (x')_b^{a(n)} ) + C5_b^n ( (\tilde{x}')_b^n - (\tilde{x}')_b^{a(n)} ) + C6_b^n (\tilde{x}')_b^n  \bigr ) \label{of-1}\\
+ &      \sum_{n\in\N}w^n d_{e^n} \sum_{q\in\Q_{e^n}}w_{q} m_{t_q}
  \Bigl (
      C^G_q (z^G)_q^n
   +  \sum_{b\in\BB} \epsilon_b ((y^+)_{b,q}^n + (y^-)_{b,q}^n) \nonumber
      \\
& \qquad
  + \Big[
\sum_{i\in\I} C^R_{i,q} (z^R)_{i,q}^n
 - P_q \bigl ( \sum_{i\in\I}( (\tilde{z}^R)_{i,q} F_i \tilde{x}_i^n - (z^R)_{i,q}^n ) \bigr )
     \Big]_{:t_q\in\T_{e^n}^R}
     \Bigr )\label{of-3} \\
- & \sum_{n\in\N_{E}}w^n \Bigl ( \sum_{i\in\I} V2_i^n \tilde{x}_i^n + \sum_{b\in\BB} V5_b^n (\tilde{x}')_b^n \Bigr ) \label{of-4},
\end{align}
\end{subequations}
\bcbl
where $w^n$ and $w_q$ denote the weights of strategic node $n$ and operational node $q$, respectively.
\bcb

The expression \eqref{of-0} gives the strategic costs associated with PV technology, which are computed as the sum of fixed preparation costs (incurred only at the first installation of a PV technology), plus variable installation costs of PV panels, plus maintenance costs.
Similarly, the term \eqref{of-1} computes the strategic costs related to BESS by considering preparation, installation and maintenance costs.
The expression \eqref{of-3} computes the expected operational costs of the domestic RES as the import costs from the power grid, plus the operating costs of BESS, plus PV electricity generation costs, minus revenue of exporting electricity to the grid.
Note that, in order to compute the total operational costs, the hourly cost of each operational node $q \in \mathcal{Q}_{e^n}, \ n \in \mathcal{N}$ is multiplied by the number of hours $m_{t_q}$ in the corresponding daily period $t_q$ and the number of days $d_{e^n}$ in the corresponding stage $e^n$.
Finally, expression \eqref{of-4} gives the expected residual value of the RES investments.


\bigskip
\noindent
\begin{subequations}\label{stra-cons}
\textbf{Strategic constraints \eqref{stra-cons}, $\forall n\in\N$}
\begin{align}
&  x_i^n\in\{0;1\}, \, \alpha_i^n\in\{0,1\}, \quad \forall i \in \I, \label{x1-2} \\
& \alpha_i^n \leq x_i^n, \quad \forall i \in \I,  \label{x1-4} \\
&  x_i^{a(n)} \leq x_i^n, \quad \forall i\in\I,  \label{x1-3} \\
& \label{x2-2}   0 \leq \tilde{x}_i^{a(n)} \leq \tilde{x}_i^n \leq \overline{\tilde{x}}_i x_i^n,                    \quad \forall i\in\I, \\
& \label{x3-2}   \sum_{i\in\I}(x_i^n - x_i^{a(n)}) \leq 1,                                            \\
& \label{x4-2}   \sum_{i\in\I}\tilde{x}_i^n \leq \overline{\overline{\tilde{x}}},                          \\
& \label{x6-2}   \underline{\tilde{x}}\alpha_i^n \leq \tilde{x}_i^n - \tilde{x}_i^{a(n)} \leq \overline{\tilde{x}}_i \alpha_i^n,  \quad \forall i\in\I, \\
& \label{b1-2}   ({x}')_b^n\in\{0,1\}, \, \beta_b^n\in\{0,1\}, \quad \forall b \in \BB, \\
& (x')_b^{a(n)} \leq (x')_b^n, \quad \forall b \in \BB, \\
& \beta_b^n \leq (x')_b^n, \quad \forall b\in\BB, \\
& \label{b2-2}   (\tilde{x}')_b^n\in\mathbb{N}, \, (\tilde{x}')_b^{a(n)} \leq (\tilde{x}')_b^n \leq \overline{(\tilde{x}')}_b (x')_b^n, \quad \forall b\in\BB, \\
& \label{b3-2}   \sum_{b\in\BB}((x')_b^n - (x')_b^{a(n)}) \leq 1, \quad \forall b\in\BB,                                  \\
& \label{b4-2}   \sum_{b\in\BB}(\tilde{x}')_b^n \leq \overline{\overline{(\tilde{x}')}}, \quad \forall b\in\BB,                 \\
& \label{b6-2}   (\underline{\tilde{x}'})\beta_b^n \leq (\tilde{x}')_b^n - (\tilde{x}')_b^{a(n)} \leq \overline{(\tilde{x}')}_b \beta_b^n, \quad \forall b\in\BB, \\
&  \label{budget-2}  \! \sum_{i\in\I} \bigl( C1_i^n (x_i^n \!- \!x_i^{a(n)}) \!+ C2_i^n (\tilde{x}_i^n \! - \! \tilde{x}_i^{a(n)}) \bigr) \!
\! + \sum_{b\in\BB}
\!
\bigl ( C4_b^n ((x')_b^n - (x')_b^{a(n)} ) + C5_b^n ((\tilde{x}')_b^n - (\tilde{x}')_b^{a(n)}) \bigr)
\!
\leq
\!
B^n. \!
\end{align}
\end{subequations}
Constraint \eqref{x1-2} defines the binary \textit{step variable} $x_i^n$ representing the installation of technology $i\in \I$ \textit{by} strategic node $n \in \N$ (i.e., wether technology $i$ is in use in node $n$) and the \textit{impulse variable} $\alpha_i^n$ representing the installation of new panels of technology $i$ at node $n$.
Constraint \eqref{x1-4} enforces the coherence between binary variables associated with PV technology $i \in \I$, stating that new PV panels of technology $i$ can be installed (i.e., $\alpha_i^n=1$) only if the technology is in use in node $n$ (i.e., $x_i^n=1$).
Constraint \eqref{x1-3} imposes that, once installed in a given strategic node, a PV technology remains in use in all successor nodes.
Constraint \eqref{x2-2} bound the number $\tilde{x}_i^n$ of PV panels of technology $i\in\I$ installed
\textit{by} strategic node $n\in\N$ from below by the number of PV panels installed \textit{by} the ancestor node $a(n)$ and from above by $\overline{\tilde{x}}_i$, which represents the maximum number of PV panels of type $i$ that can be installed in the building complex.
If PV technology $i$ is not in use in strategic node $n$ (i.e., $x_i^n=0$), constraint \eqref{x2-2} sets to zero the variable $\tilde{x}_i^n$.
Constraint \eqref{x3-2} prevents the installation of more than one PV technology at any strategic node.
Constraint \eqref{x4-2} upper bounds the total number of PV panels in use at any strategic node by $\overline{\overline{\tilde{x}}}$, representing the overall installation limit in the building complex.
Constraint \eqref{x6-2} lower and upper bounds the number of PV panels of any technology to be installed \textit{at} any strategic node, with $\underline{\tilde{x}}$ representing a minimum number of PV panels to be installed.
Constraint system \eqref{b1-2}--\eqref{b6-2} operates analogously for BESS technologies, where $\overline{(\tilde{x}^\prime)}_b$ denotes the maximum number of battery units of technology $b$ that can be installed, $\overline{\overline{(\tilde{x}^\prime)}}$ denotes the overall upper bound on battery installations, and $\underline{\tilde{x}^\prime}$ denotes a minimum number of BESS units to be installed.
Constraint \eqref{budget-2} imposes a budget $B^n$ on the RES investments at each strategic node $n$.
\bcbl
Notice that the integrality of variable $\tilde{x}_i^n$ has been relaxed due to the potentially large upper bound $\overline{\tilde{x}}_i$, whereas $(\tilde{x}')_b^n$ is kept integer.
\bcb

\bigskip
\noindent
\begin{subequations}\label{stoch-2nd-cons}
\textbf{PV and BESS operational constraints \eqref{stoch-2nd-cons}, $\forall n\in\N$}
\begin{align}
& \label{c1}    \displaystyle 0 \leq (z^R)_{i,q}^n \leq (\tilde{z}^R)_{i,q} F_i \tilde{x}_i^n,
                                           \quad \forall q\in\Q_{e^n} : t_q\in\T_{e^n}^R, \, i\in\I,\\
& \displaystyle  \Big[ \sum_{i\in\I}(z^R)_{i,q}^n \Big]_{:t_q\in\T_{e^n}^R}
                + (z^G)_q^n + \sum_{b\in\BB}((y^-)_{b,q}^n - (y^+)_{b,q}^n) = L_q + \sum_{j\in\J_1: \, t_q\in(\T^1)_{j,e^n}}
                               \bigl ( (L^1)_{j,q} - (\Delta\ell^1)_{j,q}^n \bigr )   \nonumber \\
& \displaystyle  + \displaystyle \sum_{j\in\J_2} (L^2)_{j,e^n} \! \! \sum_{q'\in\tilde{\Q}_{j,e^n,q}} \delta_{j,q'}^n,
                                           \quad \forall q\in\Q_{e^n}, \,\, \text{ where} \nonumber \\
& \label{c2} \quad \tilde{\Q}_{j,e^n,q} \, = \, \{ q'\in\Q_{e^n} \, : \, \pi_{q'}=\pi_q \, \wedge \,
t_{q'}\in(\T^2)_{j,e^n} \, \wedge t_q - (m^2)_{j,e^n,t_{q'}} +1 \leq t_{q'} \leq t_q \}, \quad \forall j\in\J_2,  \\
& \label{c5}  0 \leq m_{t_q} (y^+)_{b,q}^n \leq \rho_{b,e^n}^+ (k')_b (\tilde{x}')_b^n, \quad \forall q\in\Q_{e^n}, \, b\in\BB,\\
& \label{c6}  0 \leq m_{t_q} (y^-)_{b,q}^n \leq \rho_{b,e^n}^- (1-(f^l)_{b,e^n}) \, y_{b,a(q)}^n, \quad \forall q\in\Q_{e^n}, \, b\in\BB, \\
& \label{c7}  (1-(f^l)_{b,e^n}) \, y_{b,a(q)}^n - m_{t_q} ( (y^-)_{b,q}^n + (y^+)_{b,q}^n ) \, = \, y_{b,q}^n,
                                           \quad \forall q\in\Q_{e^n}:t_q>\underline{t}_{e^n}, \, b\in\BB,\\
& \label{c8}  0 \leq y_{b,q}^n \leq (k')_b (\tilde{x}')_b^n,
                                           \quad \forall q\in\Q_{e^n}, \, b\in\BB,
\end{align}
\end{subequations}
\bcbl
where parameter $(m^2)_{j,e^n,t_q}$ denotes the number of consecutive daily periods required to fully satisfy the electricity demand of deferrable load type $j\in\J_2$ when the supply starts at daily period $t_q\in(\T^2)_{j,e^n}$ (see Table A.4 in Appendix A for details),
given the number of hours $m_{t_q}$ in each daily period $t_q$.
The auxiliary set $\tilde{\Q}_{j,e^n,q} \subseteq \Q_{e^n}$ collects the operational nodes $q'$ such that the starting of the supply of load type $j$ in node $q'$ implies that the load must also be supplied in node $q$.
\bcb

Equation \eqref{c1} states that the share of electricity generated from PV panels to supply the RES load is bounded from above by the total electricity generated from the PV panels.
Constraint \eqref{c2} is the electricity balance equation for operational node $q$, enforcing the equality between electricity sources and uses. Specifically, the left-hand-side of \eqref{c2} computes the total electricity available in operational node $q$ as the sum of the electricity generated from the PV panels, plus the import from public power grid, plus the electricity discharged from BESS, minus the electricity charged in BESS.
The right-hand-side of \eqref{c2} gives the total electricity consumption in operational node $q$, which is the sum of the non-controlled electricity load $L_q$, plus the consumption of elastic load types $j\in\J_1$, plus the electricity requirements of deferrable load types $j\in\J_2$ supplied in operational node $q$.
In stage $e\in\E$, once a deferrable load type $j\in\J_2$ is started, it must be supplied for $(m^1)_{j,e}$ consecutive hours in order to fully satisfy its requirement.
Notice that load type $j$ must be supplied in operational node $q\in\Q_{e^n}$ if and only if the chosen starting period $t_{q'}$ falls within the interval $[t_q - (m^2)_{j,e^n,t_{q'}} + 1, \, t_q]$.
In constraint \eqref{c2}, this relationship is enforced through the set $\tilde{\Q}_{j,e^n,q}$, which includes operational nodes $q'$ belonging to the same operational scenario (i.e., $\pi_{q'}=\pi_q$), corresponding to admissible starting periods $t_{q'}\in(\T^2)_{j,e^n}$, and satisfying the above time interval condition.
Constraints\eqref{c5} and \eqref{c6} give the net bounds of the electricity charge and discharge in BESS,
by taking into account the loss and the discharge and charge depths of the battery capacity, respectively.
Constraint \eqref{c7} enforces the flow balance in the BESS for any operational node $q$
whose daily period is not the first one in the stage.
Constraint \eqref{c8} bounds the electricity storage in BESS.

\bigskip
\noindent
\begin{subequations}\label{stoch-2nd-delta}
\textbf{Load operational constraints \eqref{stoch-2nd-delta}, $\forall n\in\N$}
\begin{align}
& \big| (L^1)_{j,q}  \! - (\Delta\ell^1)_{j,q}^n  \! - \! \big( (L^1)_{j,a(q)} \! - (\Delta\ell^1)_{j,a(q)}^n \big) \big| \! \leq \overline{\Delta}r_{j,t_q} , \! \! \!  \quad \forall j\in\J_1, \,
q \in \Q_{e^n}: \, \! t_q \in (\T^1)_{j,e^n} \! \setminus \!\{\underline{t}_{e^n}\}, \label{d2} \\
& \label{d3}    0 \leq (\Delta\ell^1)_{j,q}^n \leq (\overline{\Delta}\ell^1)_{j,t_q},
                                            \quad  \forall j\in\J_1: \, t_q\in(\T^1)_{j,e^n}, \, q\in\Q_{e^n},\\
& \label{d1} \delta_{j,q}^n\in\{0,1\}, \quad \forall q\in\Q_{e^n} \, : \, t_q\in(\T^2)_{j,e^n}, \, j\in\J_2, \\
& \label{d4} \displaystyle \sum_{q\in\Q_e:
 \, t_q\in(\T^2)_{j,e^n} \wedge \pi_q = \pi'} \delta_{j,q}^n = 1,
                                      \quad \forall \pi'\in\Pi_{e^n}, \, j\in\J_2, \\
& \nonumber \delta_{j,q}^n + \delta_{j',q'}^n \leq 1, \quad \forall (j,j')\in\H_1, \, q\in\Q_{e^n}, \, q'\in\Q_{e^n}: \pi_ q=\pi_{q'} \wedge t_q\in(\T^2)_{j,e^n} \wedge t_{q'}\in(\T^2)_{j',e^n}
\\ & \quad \wedge [t_q, t_q + (m^2)_{j,e^n,t_q} - 1] \cap
                  [t_{q'}, t_{q'} + (m^2)_{j',e^n,t_{q'}} - 1] \neq \emptyset, \,  \label{d6} \\
& \nonumber \displaystyle \delta_{j,q}^n \leq \sum_{q'\in\A'_{j,j',q}} \delta_{j',q'}^n,
                        \quad \forall (j,j')\in\H_2, \, q\in\Q_{e^n}: \, t_q\in(\T^2)_{j,e^n}, \nonumber \\
& \quad \text{where } \A'_{j,j',q} \, = \, \{ q'\in\Q_e \, : \, \pi_{q'}=\pi_q \wedge t_{q'}\in(\T^2)_{j',e^n} \wedge
   t_{q'} \geq t_q + (m^2)_{j,e^n} + \overline{\Delta}t_{j,j'} \}. \label{d8}
\end{align}
\end{subequations}
Constraint \eqref{d2} is a ramp constraint for elastic loads, limiting the absolute variation of the adjusted setpoint between two consecutive daily periods.
The parameter $\overline{\Delta}r_{j,t_q}$ represents the maximum allowable difference between $(L^1)_{j,q} - (\Delta\ell^1)_{j,q}^n$ and its value in the previous daily period.
Constraint \eqref{d3} bounds the curtailment of elastic loads through the parameter $(\overline{\Delta}\ell^1)_{j,t_q}$, which denotes the maximum decrease deviation from the setpoint of load type $j$ at daily period $t_q$.
Constraint \eqref{d4} ensures that, in each operational scenario $\pi'\in\Pi_{e^n}$, exactly one daily period is selected as the starting period for deferrable load type $j\in\J_2$.
Constraint \eqref{d6} enforces incompatibility between deferrable load types.
For each pair $(j,j')\in\H_1$, it prevents overlapping activation intervals whose lengths are determined by the parameter $(m^2)_{j,e^n}$, representing the number of consecutive daily periods required to satisfy load type $j$ when started at a given period.
Constraint \eqref{d8} models minimum latency requirements between deferrable load types, by stating that if load type $j$ is supplied starting at daily period $t_q$, then load $j'$ must be satisfied at any daily period no earlier than $t_q + (m^2)_{j,e^n} + \overline{\Delta}t_{j,j'}$, with parameter $\overline{\Delta}t_{j,j'}$ denoting the required minimum separation between the completion of load $j$ and the start of load $j'$.
The auxiliary set $\A'_{j,j',q}$ collects operational nodes $q'$ in the same operational scenario that satisfy this temporal separation condition.

\bigskip
\noindent
\textbf{Tactical constraints modeling inter-stages BESS electricity flows \eqref{tac-battery_flow_balance}}
\begin{subequations}\label{tac-battery_flow_balance}
\begin{align}
 \label{tac-flow-0}
y_{b,q}^0 = &
m_{t_q} \, (y^+)_{b,q}^0 ,                             \quad \forall q\in\Q_1:t_q=1, \ b\in\BB\\
y_{b,q}^n  = &
    \frac{1}{d_{e^n}}             \sum_{q'\in\Q_{e^{a(n)}}:t_{q'}=\overline{t}_{e^{a(n)}}} w_{q'} (1-(f^l)_{b,e^{a(n)}}) \, y_{b,q'}^{a(n)}
  + \frac{d_{e^n}-1}{d_{e^n}}
  \sum_{q' \in \Q_{e^{n}}:t_{q'}=\overline{t}_{e^{n}}}          w_{q'} (1-(f^l)_{b,e^n}) \, y_{b,q'}^n \qquad \qquad \nonumber \\ &
 \label{tac-flow-3}
     + \, m_{t_q} \bigl (- (y^-)_{b,q}^n + (y^+)_{b,q}^n \bigr ),  \quad \forall q\in\Q_{e^n}:t_q = \underline{t}_{e^n}, \, n\in\N\setminus\{0\}, \, b\in\BB.
\end{align}
\end{subequations}
Equation \eqref{tac-flow-0} states that each battery $b\in\BB$ is empty at the beginning of the planning horizon.
Accordingly, the storage level in battery $b$ at the first daily period $t_q=1$ of stage $e=1$ equals the charged energy in that period.
Equation \eqref{tac-flow-3} enforces the storage electricity balance for the first daily period of any subsequent stage $e \in \mathcal{E}\setminus\{1\}$.
In order to represent the operation of long-term storage, the constraint allows for the electricity transfer between stages, by computing the storage level $y_{b,q}^{n}$ in the first daily period $t_q=\underline{t}_{e^n}$ as the sum of the following terms:
\begin{itemize}
\item $\sum_{q'\in\Q_{e^{a(n)}}:t_{q'}=\overline{t}_{e^{a(n)}}} w_{q'} (1-(f^l)_{b,e^{a(n)}}) \, y_{b,q'}^{a(n)}$ is the expected electricity level of the battery at the end of the last daily period $\overline{t}_{e^{a(n)}}$ of the previous stage $e^{a(n)}$, weighted with $\frac{1}{d_{e^n}}$ to represent the (unique) first day of stage $e^n$.
\item $\sum_{q' \in \Q_{e^{n}}:t_{q'}=\overline{t}_{e^{n}}} w_{q'} (1-(f^l)_{b,e^n}) \, y_{b,q'}^n$ is the expected electricity level of the battery at the end of the last daily period $\overline{t}_{e^{n}}$ of stage $e^{n}$, weighted with $\frac{d_{e^n}-1}{d_{e^n}}$ to represent the remaining $d_{e^n}-1$ days of stage $n$.
\item $m_{t_q} \bigl (- (y^-)_{b,q}^n + (y^+)_{b,q}^n \bigr )$
is the net battery flow in operational node $q$.
\end{itemize}
This approximation enables modeling long-term storage dynamics through the tactical multi-horizon structure, while controlling the dimensionality of the multistage stochastic tree. (see Fig. \ref{MH_Tree_b}).

\bigskip
\noindent
\textbf{Maximum expected discomfort bound \eqref{rn},
$\forall n\in\N$}
\begin{align}
\label{rn} \displaystyle
\sum_{\pi\in\Pi_{e^n}}
w'_\pi
\sum_{q\in\A_\pi}
    \Bigl (
    \sum_{j\in\J_1: \,   t_q\in(\T^1)_{j,e^n}}
    m_{t_q}
    (D^1)_{j,t_q} (\Delta\ell^1)_{j,q}^n
  + \sum_{j\in\J_2: \, t_q\in(\T^2)_{j,e^n}} (D^2)_{j,t_q} \, \delta_{j,q}^n
    \Bigr )
\leq \hat{D}_{e^n}
,
\end{align}
where $w'_\pi$ denotes the probability of operational scenario $\pi$.

To upper bound user discomfort arising from deviations in household electricity consumption from the elastic setpoints, as well as shifts in the supply of deferrable loads relative to the reference period, constraint \eqref{rn} establishes an upper bound $\hat{D}_{e^n}$ on the expected discomfort at each strategic node.
It is referred to as a \textit{risk-neutral} approach to discomfort.

While the risk-neutral approach bounds the expected discomfort in each strategic node, it does not prevent rare but severe operational outcomes, leaving the system exposed to extreme user discomfort.
To explicitly control such low-probability, high-impact events, we adopt first- and second-order SD measures, which have been extensively studied in the stochastic optimization literature as coherent and time-consistent criteria (see, e.g., \cite{dentchevaruszczynski03,Gollmer08,Gollmer11,forestry18,dismit18}).
These measures allow direct control of both the probability and the magnitude of discomfort exceedances beyond predefined thresholds.

{\color{black}
Before introducing the SD constraint system, we briefly relate this choice to alternative risk measures commonly used in stochastic optimization. Value-at-Risk (VaR) controls a quantile of the cost (or discomfort) distribution but does not restrict the magnitude of realizations beyond the selected confidence level, potentially allowing very large costs in low-probability scenarios. Conditional Value-at-Risk (CVaR) partially addresses this issue by accounting for expected tail costs; however, it requires a careful calibration of weighting parameters in the objective function, which may reduce transparency in the resulting trade-offs. In contrast, stochastic dominance (SD) constraints provide a more explicit and modeler-driven control of risk by directly bounding both the probability and the severity of threshold exceedances. This structure makes SD-based formulations particularly suitable for controlling user discomfort in domestic RES planning, where rare but extreme operational scenarios must be limited independently of objective-function trade-offs.
SD constraints are also closely tied to utility theory, see \cite{dentchevaruszczynski03}.
This connection is particularly valuable in buildings with multiple residents and commercial occupants, each with distinct utility preferences. SD constraints naturally incorporate these heterogeneous preferences in a unified framework, providing a more appropriate approach than relying on a single VaR or CVaR measure calibrated for one central decision maker.
}

\bigskip
\noindent
\textbf{Risk-averse discomfort constraint system \eqref{sd},
$\forall p\in\P_{e^n}, \, n\in\N$}
\begin{subequations}
\label{sd}
\begin{align}
&  \sum_{q\in\A_\pi}
    \Bigl ( \! \!
    \sum_{j\in\J_1: \, t_q\in(\T^1)_{j,e^n}}
   \! \! \! \!
    m_{t_q}
    (D^1)_{j,t_q} (\Delta\ell^1)_{j,q}^n
  +
  \! \! \! \!
  \sum_{j\in\J_2: \, t_q\in(\T^2)_{j,e^n}}
  \! \! \! \!
  (D^2)_{j,t_q} \delta_{j,q}^n
    \Bigr )
  - s^{p,\pi} \leq \overline{D}^p, \, \forall \pi\in\Pi_{e^n},
  \label{sd-1}
\\
& \label{sd-2}  0 \leq s^{p,\pi} \leq \overline{s}^p \overline{D}^p \eta^{p,\pi},
\quad \forall \pi\in\Pi_{e^n},
\\
& \label{sd-3}  \displaystyle
   \eta^{p,\pi}\in\{0,1\}, \quad \forall \pi\in\Pi_{e^n}, \\
& \label{sd_3bis}
   \sum_{\pi\in\Pi_{e^n}} w'_\pi \eta^{p,\pi} \leq \overline{\eta}^p 
   ,
\\
& \label{sd-4}  \displaystyle
   \sum_{\pi\in\Pi_{e^n}} w'_\pi s^{p,\pi} \leq \overline{\overline{s}}^p \overline{D}^p.
\end{align}
\end{subequations}
Constraint \eqref{sd-1} defines, for each operational scenario $\pi \in \Pi_{e^n}$, the decision variable $s^{p,\pi}$ as the daily discomfort excess beyond the threshold $\overline{D}^p$.
Constraint \eqref{sd-2} sets an upper bound to the discomfort excess and ensures consistency between the variables $s^{p,\pi}$ and $\eta^{p,\pi}$, which indicates whether the discomfort threshold is violated in scenario $\pi$.
Constraint \eqref{sd-3} defines the domain of the variables $\eta^{p,\pi}$.
Constraint \eqref{sd_3bis} bounds the probability of exceeding the threshold,
while constraint \eqref{sd-4} bounds the expected exceedance in each strategic node.
Specifically, the first-order SD constraints are given by constraint system \eqref{sd-1}--\eqref{sd_3bis},
while the second-order SD measure is given by \eqref{sd-1}--\eqref{sd-2}, \eqref{sd-4}.
It is worth pointing out that the first-order SD measure has a direct relationship with the time-honored chance-constrained approach introduced in the seminal work of \cite{Charnes58}, where the probability of violating a profile threshold is upper bounded.

\subsection{Compact formulation}\noindent
For compactness purposes, let $\textbf{x}^n$ denote the strategic variables in node $n\in\N$ and $\textbf{y}^n_q$ denote the operational variables in node $q\in\Q_{e^n}, \, n\in\N$, where
%
\begin{align*}
 \textbf{x}^n = & \big[ x_i^n, \, \tilde{x}_i^n, \, \alpha_i^n, \,  (x')_b^n, \,  (\tilde{x}')_b^n, \, \beta_b^n \big], \, \text{with } i\in\I, \, b\in\BB, \\
 \textbf{y}^n_q = &
\Bigl [
\big( (z^R)_{i,q}^n \big) _{:t_q\in\T_{e^n}^R},  \, (z^G)_q^n, \, (\Delta\ell^1)_{j,q}^n, \,
\big( \delta_{j',q}^n \big)_{:t_q\in(\T^2)_{j,e^n}},
(y^-)_{b,q}^n, \, (y^+)_{b,q}^n, \, y_{b,q}^n, \,
\big( y_{b,q'}^{a(n)} \big )_{q'\in\Q_{e^{a(n)}} \, \wedge \, t_{q'}=\overline{t}_{e^{a(n)}}}
\Big] \, \\
& \text{with }
 i\in\I, \, b\in\BB, \, j\in\J_1, \, j'\in\J_2.
\end{align*}

Additionally, let the vector $\textbf{s}^{p,\pi}$ denote the operational scenario variables for risk-averse discomfort modeling in operational scenario $\pi\in\Pi_{e^n}$ for policy profile $p\in\P_{e^n}$ in strategic node $n\in\N$:
\begin{align*}
\textbf{s}^{p,\pi} = [ s^{p,\pi}, \, \eta^{p,\pi}].
\end{align*}
Let $a^n$ and $b^n_q$ denote the vectors of the objective function coefficients
for the elements of the vectors of variables $\textbf{x}^n$ and $\textbf{y}^n_q$, respectively.
The RES model \eqref{of}$-$\eqref{sd} can be expressed in a compact form as
\begin{subequations}\label{syn}
\begin{align}
\label{syn-1}
z^* \, = \, \min  \, & \sum_{n\in\N} w^n \Big(a^n \textbf{x}^n +
\sum_{q\in\Q_{e^n}} w_q \,b^n_q  \,\textbf{y}^n_q \Big) \\
\label{syn-2} \text{s.t. } & \text{Strategic constraints\eqref{stra-cons} in } (\textbf{x}^{a(n)}, \, \textbf{x}^n), \qquad & \forall n\in\N,\\
& \label{syn-3}
\text{PV and BESS operational constraints \eqref{stoch-2nd-cons} in } (\mathbf{x}^n, \, \textbf{y}^n_q), \qquad &\forall  q\in\Q_{e^n}, \, n\in\N, \\
& \label{syn-4}
\text{Load operational constraints \eqref{stoch-2nd-delta} in } \mathbf{y}^n_q, \qquad &\forall  q\in\Q_{e^n}, \, n\in\N, \\
& \label{syn-5}
\text{Tactical constraints \eqref{tac-battery_flow_balance} in } \mathbf{y}^n_q, \qquad &\forall  q\in\Q_{e^n}, \, n\in\N, \\
& \label{syn-6ub}
\text{Maximum expected discomfort constraints \eqref{rn} in } \mathbf{y}^n_q, \qquad &\forall  n\in\N, \\
& \label{syn-6}
\text{Discomfort SD constraints \eqref{sd} in } \big( (\textbf{y}^n_q)_{:\pi_q=\pi},\textbf{s}^{p,\pi}
\big), \qquad &
\forall \pi\in\Pi_{e^n}, \, n\in\N.
\end{align}	
\end{subequations}

\section{A matheuristic algorithm for obtaining feasible solutions}\label{sec:HEUR}

In realistic applications, the number of strategic and operational nodes in the RES model \eqref{of}--\eqref{sd} grows rapidly with the number of stages and scenarios.
As a consequence, solving large instances of the multistage stochastic MILP to proven optimality using a state-of-the-art solver becomes computationally impractical.
To address this challenge, we develop a structured matheuristic approach, denoted SFR3 (Scenario variables Fixing and Randomization of constraints / variables' integrality Relaxation iterative Reduction).
The algorithm builds upon the decomposition principles introduced in \cite{shubnep21}, but is here specifically adapted to the multi-horizon RES framework.
The key idea of SFR3 is to approximate the original multistage stochastic program by solving a sequence of reduced multi-horizon subproblems.
At each iteration, only a subset of strategic nodes is considered.
Future stages may be partially relaxed through probabilistic node selection.
The decisions corresponding to the earliest non-relaxed stage are then fixed.
The horizon is shifted forward and the process is repeated until all stages have been processed.
\textcolor{black}{
By construction, the SFR3 procedure guarantees feasibility for the original RES problem.
At each iteration, a reduced subproblem is solved and the decisions associated with the earliest non-relaxed stage are fixed.
By iterating this process across all stages, a complete feasible solution to the original problem is obtained.
By fixing decisions stage by stage through consistent submodel solutions, a complete feasible solution is obtained.
Due to selective node inclusion and probabilistic relaxation, global optimality cannot be guaranteed.
The quality of the solutions produced by SFR3 is therefore assessed through lower bounds and benchmarking against alternative decomposition heuristics, as illustrated in the next sections.
}

\textcolor{black}
{
This section is organized as follows.
Section \ref{sec:Algo} describes in detail the structure and main steps of the SFR3 algorithm, while Section \ref{sec:choice_SFR3} discusses the parameter strategies and provides guidance for selecting the input parameters.
}

\subsection{Algorithmic structure}
\label{sec:Algo}
Algorithm \ref{algo} provides a detailed description of the procedure designed to compute a feasible solution to the RES problem, with the corresponding notation reported in Table \ref{tab:SFR3_par}.

\begin{longtable}[ht!] {p{2.2 cm} p{13.0 cm}}
    \hline
    Symbol & Description \\
    \hline
    $\kappa$   &  iteration number of the algorithm. \\

    $\hat{e}$  &  number of non-relaxed stages, which nodes are not subject to any relaxation at the iterations.
     So, $\kappa=\{1,\ldots,E-\hat{e}+1)$. \\

    $\hat{e}^R$ & number of relaxation stages, which nodes can be randomly selected at each iteration. \\


   $\phi_{e}$ &  probability of selecting the strategic nodes belonging to the relaxation stage $e$. \\

   $\S^n\subseteq\N$ & set of successors of strategic node $n\in\N$. \\

  $\S^n_1\subseteq\S^n$ & set of children nodes of strategic node $n$. \\

   $\mathbf{x}^n$ & vector of the strategic decision variables defined in strategic node $n$. \\

   $\mathbf{y}^n_q$ & vector of the operational decision variables defined in operational node $q\in\Q_{e^n}$. \\

   $\mathbf{s}^{p,\pi}$ & vector of decision variables for discomfort modeling in the operational scenario $\pi\in\Pi_{e^n}$
   for policy profile $p\in\P_{e^n}$. \\

   $\overline{\S}^n\subseteq\S^n$ & set of successors of strategic node $n$ that belong to a non-relaxed stage at an iteration. \\

   $\tilde{\N}_{e}^n \subseteq\S^n \cap\N_{e}$ & set of successors of strategic node $n$ that belong to relaxation stage $e$.
   Those nodes are randomly selected to be included in the model at an iteration. \\

   $\overline{\overline{\S}}^n \subseteq \S^n$ & set of successors of strategic node $n$ that belong to any relaxation stage.
   Those nodes are randomly selected to be included in the model at an iteration. \\

   $\overline{\N}^\kappa \subseteq \N$ & set of all strategic nodes to be considered at iteration $\kappa$. \\

   $\hat{w}^n$ & updated weight of node  $n\in\overline{\N}^{\kappa}$ at iteration $\kappa$. \\
\hline
      \caption{Parameters and sets in the SFR3 approach}
    \label{tab:SFR3_par}
\end{longtable}
%
%
\begin{algorithm}[ht]
\caption{SFR3 algorithm to compute a feasible solution to the RES problem}
\begin{algorithmic}[1]
\label{algo}
\STATE Input: $\hat{e}, \, \hat{e}^R, \,\phi_{e}$
\FOR{$\kappa = 1$ to $E-\hat{e}+1$}
   \FOR{$r\in\N_\kappa$}
     \STATE Construct set $\overline{\S}^r:=\{n\in\S^r:e^n<\kappa+\hat{e}\}$
     \IF{$\hat{e}^R>0$}
     \FOR{$e'=\kappa+\hat{e}$ to $\min(\kappa+\hat{e}+\hat{e}^R-1,E)$}
        \STATE Construct set $\tilde{\N}^r_{e'}$ by randomly selecting with probability $\phi_{e'}$  strategic nodes $n'\in S^r \cap \N_{e'}$,
        with $a(n') \in \overline{\S}^r \cup \tilde{\N}^r_{e'-1}$
     \ENDFOR
     \vspace{0.1cm}
     \STATE  $\overline{\overline{\S}}^r := \bigcup_{e'=\kappa}^{\min(\kappa+\hat{e}+\hat{e}^R-1,E)} \tilde{\N}^r_{e'}$
    \ELSE
    \vspace{0.1cm}
    \STATE $\overline{\overline{\S}}^r := \emptyset$
    \ENDIF
     \STATE Define set $\overline{\N}^r := \{r\} \cup \overline{S}^r \cup \overline{\overline{\S}}^r$
     \STATE Set $\hat{w}^r=1$
     \FOR{$n\in {\overline{\N}}^r \setminus\{r\}$}
     \STATE $\hat{w}^n =\frac{\hat{w}^{a(n)}}{\sum_{m\in \S_1^{a(n)} \cap \tilde{\N}^r_{e^n}  }w^m} w^n$
     \ENDFOR
     \STATE Solve the $r$-submodel, defined as a replica of model \eqref{of}$-$\eqref{sd}, restricted to set ${\overline{\N}}^r$, with weights  $\hat{w}^n,\, n \in{\overline{\N}}^r$, and obtain the solution of strategic variables $\textbf{x}^n$, operational variables $\textbf{y}^n_q, \ q\in\Q_{e^n}$ and discomfort variables $\mathbf{s}^{p,\pi}, \, \pi\in\Pi_{e^n}, \, p\in\P_{e^n}$, for $n\in\overline{\N}^r$
     \STATE Fix variables $\textbf{x}^r$, $\textbf{y}^r_q$, $q\in\Q_{e^r}$, and $\mathbf{s}^{p,\pi}$, $\pi\in\Pi_{e^r}, \, p\in\P_{e^r}$,
     to the solution of the $r$-submodel
     \IF{$\kappa=E-\hat{e}+1$}
     \STATE Fix variables $\textbf{x}^n$ and $\textbf{y}^n_q$, $q\in\Q_{e^n}$, and $\mathbf{s}^{p,\pi}$, $\pi\in\Pi_{e^n}, \, p\in\P_{e^n}$, for $n\in\overline{\S}^r$ to the solution of the $r$-submodel
     \ENDIF
\ENDFOR
\ENDFOR
\STATE Output: A feasible solution
$( \mathbf{x}^n, \, \mathbf{y}^n_q \, \forall q\in\Q_{e^n}, \, \mathbf{s}^{p,\pi} \, \forall \pi\in\Pi_{e^n}, \, p\in\P_{e^n} )$ $\forall n\in\N$.
\end{algorithmic}
\end{algorithm}
Specifically, at a given iteration $\kappa$ of the algorithm, after fixing the decisions of stages $1,\ldots,\kappa-1$, the set $\E$ of stages is partitioned into the following three disjoint \textit{stage groups}:
\begin{itemize}
    \item \textit{Non-relaxed} stages $\{\kappa,\ldots,\kappa+\hat{e}-1\}$, whose strategic and operational nodes are fully included in the subproblem.
    \item \textit{Relaxation} stages $\{\kappa+\hat{e},\ldots,\kappa+\hat{e}+\hat{e}^R-1\}$, where only a subset of strategic nodes is selected according to probability parameter $\phi_e$.
    \item \textit{Removed} stages $\{\kappa+\hat{e}+\hat{e}^R,\ldots,E\}$, which are excluded from the subproblem at iteration $\kappa$.
\end{itemize}
Having fixed the values of decision variables up to stage $\kappa-1$,
the MH problem is decomposed into $|\N_\kappa|$ independent subproblems, each rooted at a strategic node $r \in \mathcal{N}_\kappa$.
Each $r$-submodel is a restricted replica of the RES model \eqref{of}--\eqref{sd}, defined over a reduced set of strategic nodes.
For each root node $r \in \mathcal{N}_\kappa$, the corresponding $r$-submodel is constructed by identifying the strategic nodes to be included according to the stage partition defined above.
The following sets are defined to specify this reduced strategic tree.
\begin{itemize}
    \item Successors in non-relaxed stages:
    $\overline{\mathcal{S}}^r = \{ n \in \mathcal{S}^r : e^n < \kappa + \hat{e} \}$.
    These nodes belong to the non-relaxed stages and are fully included in the submodel.

    \item Selected successors in relaxation stages:
    $\tilde{\mathcal{N}}_{e'}^r \subseteq \mathcal{S}^r \cap \mathcal{N}_{e'}$ for stages
    $\kappa + \hat{e} \leq e' \leq \min(\kappa+\hat{e}+\hat{e}^R-1, E)$.
    These nodes are randomly selected with probability $\phi_{e'}$, subject to path consistency, i.e., a node is selected only if its ancestor is also selected.

    \item All selected nodes in relaxation stages:
    $\overline{\overline{\mathcal{S}}}^r = \bigcup_{e'} \tilde{\mathcal{N}}_{e'}^r$.
    This set aggregates all strategic nodes selected in the relaxation stages.
    If $\hat{e}^R = 0$, this set is empty.

    \item Complete set of strategic nodes in the submodel:
    $\overline{\mathcal{N}}^r = \{ r \} \cup \overline{\mathcal{S}}^r \cup \overline{\overline{\mathcal{S}}}^r$.
    This defines the reduced strategic tree on which the $r$-submodel is solved.
\end{itemize}

To preserve probability consistency in the reduced scenario tree, node weights are rescaled.
The root node $r$ of the $r$-submodel is assigned weight $\hat{w}^r=1$.
The weight of each selected successor node $n\in\overline{\N}^r \setminus\{r\}$ is updated recursively so that each node's weight equals the sum of the weights of its children in the reduced tree, as computed in line 16 of Algorithm \ref{algo}.
This guarantees probability consistency in the reduced scenario tree despite node reduction.

The resulting $r$-submodel is solved to determine the strategic variables $\mathbf{x}^n$, operational variables $\mathbf{y}^n_q$, and discomfort variables $\mathbf{s}^{p,\pi}$ for all nodes included in $\overline{\mathcal{N}}^r$.
The decisions corresponding to stage $\kappa$ are then fixed in the original problem.
The algorithm proceeds to iteration $\kappa+1$ until all stages have been processed.

\subsection{Parameter strategies}
\label{sec:choice_SFR3}
{\color{black}
The performance of SFR3 depends on the parameters $\hat{e}$, $\hat{e}^R$, and $\phi_e$.
The parameter $\hat{e}$ determines the number of non-relaxed stages included in each subproblem and therefore controls the look-ahead horizon.
The parameter $\hat{e}^R$ determines the number of relaxation stages in which only a subset of strategic nodes is included.
The parameter $\phi_e$ controls the probability of selecting nodes in relaxation stages and regulates the trade-off between scenario representativeness and computational effort.
Different parameter configurations correspond to progressively richer horizon approximations.
The main strategy classes are summarized below.
\begin{itemize}
    \item The \textit{weak myopic strategy}: $\hat{e}=1$, $\hat{e}^R=0$, $\phi_e=0$.
    Each subproblem only includes the current stage.
    This minimizes computational effort but ignores inter-stage dependencies and may lead to suboptimal long-term decisions.

    \item The \textit{stronger myopic strategy}: $\hat{e}=2$, $\hat{e}^R=0$, $\phi_e=0$.
    The look-ahead window is extended to one additional stage.
    This improves solution quality by capturing short-term inter-stage effects, at moderate computational cost.

    \item The \textit{multistage myopic strategy}: $\hat{e}>2$, $\hat{e}^R=0$, $\phi_e=0$.
    Several future stages are fully included in each subproblem.
    This captures longer-term interactions across stages but increases subproblem size and computational time.

    \item \textit{Relaxed-horizon strategies:} $\hat{e}^R>0$.
    Future stages beyond the non-relaxed window are partially explored through probabilistic node selection governed by $\phi_e$.
    These strategies increase scenario diversity without fully expanding the subproblem dimension.
    Larger values of $\phi_e$ enhance representativeness of uncertainty at the cost of higher computational effort.
\end{itemize}
Overall, increasing $\hat{e}$ improves horizon consistency, while introducing relaxation stages through $\hat{e}^R$ and $\phi_e$ provides a controlled mechanism to balance solution quality and tractability.
}

\section{Lower bounds for multistage multi-horizon models}\label{sec:lb}

{\color{black}
Assessing the quality of feasible solutions is a critical aspect of solving large-scale stochastic optimization problems, especially when the underlying mixed-integer formulation cannot be solved to proven optimality. In such cases, lower bounds provide a valuable reference for evaluating the performance of heuristic and matheuristic solution approaches, allowing one to quantify how far a feasible solution cost may lie from the best achievable objective function value.

In this work, we consider several lower-bound formulations for the RES model \eqref{of}--\eqref{sd}. These bounds are not intended as alternative solution algorithms, but rather as benchmarking tools to assess the quality of the solutions produced by the proposed rolling-horizon matheuristic. The considered approaches differ in how uncertainty is approximated and how non-anticipativity constraints are relaxed across strategic stages and scenarios.
}
Specifically, we consider the following approaches:
\begin{enumerate}
\item $MHEV$, \textit{Multi Horizon Expected Value}, where all strategic and operational uncertain parameters are replaced by their expected values.
It is specially useful for very large-sized instances.

\item $MHOEV$, \textit{Multi Horizon Operational Expected Value},
where operational uncertain parameters are replaced by their expected values.
So, the uncertainty is only considered in the strategic stochastic parameters.

\item $SWS$, \textit{Strategic Wait-and-See}, where the non-anticipativity constraints on the strategic decision variables are relaxed.
So, the strategic scenarios are considered in an independent mode.


\item $SMG$, \textit{Scenario Multistage Grouping}, where the strategic scenarios are grouped in disjoint groups,
   so that the non-anticipativity constraints that link one scenario group with any other are relaxed.

\item $SMC$, \textit{Scenario Multistage Clustering}, where the non-anticipativity constraints are relaxed up to a given stage, so as to obtain clusters of strategic scenarios.
\end{enumerate}

{\color{black}
Lower bounds based on the $MHEV$, $MHOEV$, and $SWS$ formulations were originally introduced in \cite{Maggioni20}. For completeness, the corresponding models are reported in Appendix B, where we also provide a brief summary of their main characteristics.
In the main text, we focus on the $SMG$ and $SMC$ approaches, which represent methodological extensions beyond the existing literature of bounding techniques for multi-horizon models. These methods relax non-anticipativity constraints in a structured manner, either by grouping strategic scenarios or by clustering them up to a given stage, and provide tighter lower bounds for large-size instances.
}

\subsection{Scenario Multistage Grouping (SMG)}\label{sec:SMG}

The lower bound provided by  $SWS$ can be improved by grouping the strategic scenarios into disjoint groups.
The grouping strategy is a widely employed technique for obtaining lower bounds in stochastic problems. For its application to multistage stochastic programs, we refer the reader to \cite{Maggioni14,Maggioni16, Maggioni16CMS,Maggioni19}.
For the MH case, let $SMG$ be the method that groups the scenarios of the original problem into $G$ independent subproblems.
A smaller $G$ results in a stronger lower bound but comes at the cost of increased computational effort.
Table \ref{tab:SMG_par} reports the notation needed to formulate $SMG$.

For any scenario group $g\in\G$,
the elements of set $\Omega_g$ can be randomly chosen from set $\Omega$, without repetition (i.e., $\Omega_g\cap\Omega_{g'} = \emptyset, \, \, \forall g,g'\in\G: g \neq g'$), such that $\Omega = \bigcup_{g\in\G}\Omega_g$.
For any $g\in\G$, we define as $g$-submodel the replica of the RES model \eqref{of}--\eqref{sd} restricted to the set $\bar{\N}_g\subseteq\N$ of strategic nodes associated with strategic scenarios $\omega \in \Omega_g$.
In any $g$-submodel, the weights must be properly rescaled.
Let $\tilde{w}_\omega$ denote the updated weight of scenario $\omega$ in the corresponding group $g$, and $\tilde{w}_g^n$ does it for the  updated weight in group $g$ of a given node $n\in\bar{\N}_g$; see the computations in Table \ref{tab:SMG_par}.

\bigskip
\begin{longtable}[ht!] {p{3.8 cm} p{11.2 cm}}
    \hline
    Symbol & Description \\
    \hline
    $g\in\G=\{1,\ldots,G\}$  &  set of groups of strategic scenarios, where $G=|\G|$.  \\

    $\Omega_g \subseteq \Omega$ & set of scenarios in group $g\in\G$. \\

    $\bar{\N}_g = \bigcup_{\omega\in\Omega_g}\A^\omega$ & set of nodes in group $g\in\G$. \\

    $W_g=\sum_{\omega\in\Omega_g}w^\omega$ & weight of group $g\in\G$. \\

    $\tilde{w}_\omega = \frac{w^\omega}{\sum_{w'\in\Omega_g}w^{\omega'}}$ & weight of scenario $\omega$ in group $g$, for $\omega\in\Omega_g$.  \\
    $\tilde{w}_g^n = \sum_{\omega\in\Omega^n\cap\Omega_g} \tilde{w}_\omega$ & weight assigned to node $n \in\bar{\N}_g$ in group $g\in\G$. \\

    $\textbf{x}_g^n$ & replica of the vector of strategic decision variables $\textbf{x}^n$ in node $n$ in group $g$, for $n\in\bar{\N}_g, \, g\in\G$. \\

    $\textbf{y}_{q,g}^n$ & replica of the vector of operational decision variables $\textbf{y}^n_{q}$ in operational node $q$ in group $g$, for $q \in \Q_{e^n}, \ n\in\bar{\N}_g, \, g\in\G$. \\

    $\textbf{s}_g^{p,\pi}$ & replica of the vector of variables for discomfort excess $\textbf{s}_g^{p,\pi}$ in operational scenario $\pi$ in group $g$, for profile $p$, for  $\pi \in \Pi_{e^n}, \ p \in \P_{e^n}, \ n\in\bar{\N}_g, \, g\in\G$. \\

\hline
      \caption{Sets and parameters in $SMG$}
    \label{tab:SMG_par}
\end{longtable}
Thus, the $g$-submodel can be expressed
\begin{align*} \label{g-submodel}
z_g (G) \, = \, \min \, & \sum_{n\in\bar{\N}_g} \tilde{w}_g^n
\Big( a^n \textbf{x}_g^n + \sum_{q\in\Q_{e^n}} w_q \, b^n_q \, \textbf{y}_{q,g}^n \Big)\\
\text{s.t. } & \eqref{syn-2}, \, \forall n\in\bar{\N}_g,
\text{ replacing } (\mathbf{x}^{a(n)}, \mathbf{x}^{n}) \text{ with } (\mathbf{x}_{g}^{a(n)}, \mathbf{x}_{g}^n),\\
& \eqref{syn-3}, \,
\forall q\in\Q_{e^n}, \ n\in\bar{\N}_g,
\text{ replacing } (\mathbf{x}^{n}, \mathbf{y}^{n}_q) \text{ with } (\mathbf{x}^n_g, \mathbf{y}_{q,g}^n), \\
& \eqref{syn-4}-\eqref{syn-5} \,
\forall q\in\Q_{e^n}, \ n\in\bar{\N}_g,
\text{ replacing } \mathbf{y}^{n}_q \text{ with } \mathbf{y}_{q,g}^n, \\
& \eqref{syn-6ub}, \, \forall n\in\bar{\N}_g,
\text{ replacing } \mathbf{y}^{n}_q \text{ with } \mathbf{y}_{q,g}^n, \\
& \eqref{syn-6}, \, \forall  \pi \in \Pi_{e^n}, \, p \in \P_{e^n}, \, n\in\bar{\N}_g,
 \text{ replacing } (\mathbf{y}^{n}_q, \mathbf{s}^{p,\pi}) \text{ with } (\mathbf{y}_{q,g}^n, \mathbf{s}^{p,\pi}_g).
\end{align*}
By combining the optimal values of the scenario groups with weights $W_g=\sum_{\omega \in \Omega_g} w^\omega$ we can compute a lower bound $z_{SMG}(G)=\sum_{g\in\G} W_g \,  z_g(G)$ of optimal solution cost value $z^*$.

\subsection{Scenario Multistage Clustering (SMC)}\label{sec:SMC}
\begin{figure}[ht!]
\begin{center}
\includegraphics[height=5cm, viewport=112 522 300 675,clip=]{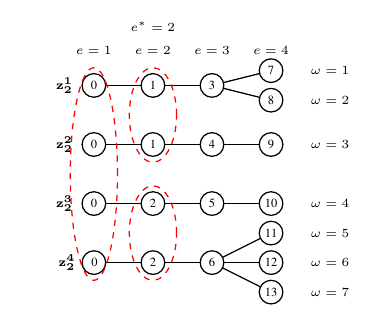}
\caption{Breaking-stage $e^*=2$ and scenario cluster sets.
$\N_1=\{0\}, \, \S^0_1=\{1,2\}, \, \N_{e^*+1}=\{3,4,5,6\}, \, n_1=3, n_2=4, n_3=5, n_4=6$}
\label{breaking-stage}
\end{center}
\end{figure}

Another strategy to strengthen the lower bound provided by $SWS$ consists of relaxing the non-anticipativity constraints for the strategic variables up to a given stage $e^* <E$, which is referred to as \textit{breaking stage}.
So, the MH problem can be decomposed into $|\mathcal{N}_{e^*+1}|$ independent submodels, referred to as the \( c \)-submodels, for \( c \in \mathcal{C} \), where \( \mathcal{C} \) represents the set of strategic scenario clusters. The cardinality of $\mathcal{C}$ equals $|\mathcal{N}_{e^*+1}|$.
It is worth mentioning that each node in $\N_{e^*+1}$ and its successors belong to just one scenario cluster in set $\C$.
$SMC$ is based on the solution of the $c$- submodels.
For its application to multistage stochastic programs, we refer the reader to \cite{mcld17}.
Table \ref{tab:SMC_par} shows the notation needed to formulate it.
It can be observed in the table that $\Omega_c\cap\Omega_{c'} = \emptyset, \, \, \forall c,c'\in\C: c \neq c'$ and
$\Omega = \bigcup_{c\in\C}\Omega_c$.
Figure \ref{breaking-stage} provides a visual representation of $SMC$ for a scenario tree consisting of 13 nodes and 7 scenarios.
By selecting the second stage as the breaking one $e^*=2$, the scenario tree decomposes into $|\N_3|=4$ scenario clusters, with $\Omega_1 = \{1,2\}$, $\Omega_2 = \{3\}$, $\Omega_3 = \{4\}$, and $\Omega_4 = \{5,6,7\}$.
%
\begin{longtable}[ht!] {p{4.4 cm} p{10.6 cm}}
    \hline
    Symbol & Description \\
    \hline
    $e^* \in \E\setminus\{E\}$ & breaking stage \\

    $c\in\C=\{1,\ldots,|\N_{e^*+1}|\}$  &  set of clusters of strategic scenarios. \\

    $n_c\in \N_{e^*+1}$  &  strategic node in stage $e^*+1$ associated with cluster $c\in\C$. \\


    $\Omega_c \subseteq \Omega$ & set of scenarios in cluster $c\in\C$. \\

    $\N'_c = \bigcup_{\omega\in\Omega_c}\A^\omega$ & set of nodes in cluster $c\in\C$. \\

    $\bar{W}_c=\sum_{\omega\in\Omega_c}w^\omega$ & weight of cluster $c\in\C$. \\

    $\bar{w}_\omega = \frac{w^\omega}{\sum_{w'\in\Omega_c}w^{\omega'}}$ & weight of scenario $\omega$ in cluster $c$,
    for $\omega\in\Omega_c$.  \\

    $\bar{w}_c^n = \sum_{\omega\in\Omega^n\cap\Omega_c} \tilde{w}_\omega$ & weight assigned to node $n \in\N'_c$ in cluster $c\in\C$. \\

    $\textbf{x}_c^n$ & replica of the vector of strategic decision variables $\textbf{x}^n$ in node $n$ in scenario cluster $c$,
    for $n\in\N'_c, \, c\in\C$. \\

    $\textbf{y}_{q,c}^n$ & replica of the vector of operational decision variables $\textbf{y}^n_{q}$ in node $n$, operational node $q$ in scenario cluster $c$, for $q \in \Q_{e^n}, \ n\in\N'_c, \, c\in\C$. \\

    $\textbf{s}^{p,\pi}_c$ & replica of the vector of variables for discomfort excess $\textbf{s}^{p,\pi}_{c}$ in operational scenario $\pi$ in cluster $c$ for profile $p$, for $\pi \in \Pi_{e^n}, \, p \in \P_{e^n}, \ n\in\N'_c, \, c\in\C$. \\
\hline
      \caption{Sets and parameters in $SMC$}
    \label{tab:SMC_par}
\end{longtable}

For any $c\in\C$, we define as $c$-submodel the replica of the RES model \eqref{of}--\eqref{sd} restricted to the set ${\N}'_c\subseteq\N$ of strategic nodes associated with strategic scenarios $\omega \in \Omega_c$.
Similar to the $SMG$ case, the weights in $c$-submodel must be properly rescaled.
Let $\bar{w}_\omega$ denote the updated weight of scenario $\omega$ in the corresponding cluster $c$, and $\bar{w}_c^n$ does it for the updated weight of  node $n\in\N'_c$ in cluster $c$; see the computations in Table \ref{tab:SMC_par}.
Thus, the $c$-submodel can be expressed
\begin{align*}
z_c (e^*)  = \min \,& \sum_{n\in{\N}'_c} \bar{w}_c^n
\Big(
a^n \textbf{x}_c^n + \sum_{q\in\Q_{e^n}} w_q \, b^n_q \, \textbf{y}_{q,c}^n
\Big)\\
\text{s.t. } & \eqref{syn-2}, \, \forall n\in{\N}'_c,
\text{ replacing } (\mathbf{x}^{a(n)}, \mathbf{x}^{n})  \text{ with } (\mathbf{x}_{c}^{a(n)}, \mathbf{x}_{c}^n),\\
& \eqref{syn-3}, \,
\forall q\in\Q_{e^n}, \ n\in{\N}'_c,
\text{ replacing } (\mathbf{x}^{n}, \mathbf{y}^{n}_q) \text{ with } (\mathbf{x}^n_c, \mathbf{y}_{q,c}^n), \\
& \eqref{syn-4}-\eqref{syn-5}, \,
\forall q\in\Q_{e^n}, \ n\in{\N}'_c,
\text{ replacing } \mathbf{y}^{n}_q \text{ with } \mathbf{y}_{q,c}^n, \\
& \eqref{syn-6ub}, \,
\forall n\in{\N}'_c,
\text{ replacing } \mathbf{y}^{n}_q \text{ with } \mathbf{y}_{q,c}^n, \\
& \eqref{syn-6}, \, \forall  \pi \in \Pi_{e^n}, \, p \in \P_{e^n}, \, n\in{\N}'_c,
\text{ replacing } (\mathbf{y}^{n}_q, \mathbf{s}^{p,\pi}) \text{ with } (\mathbf{y}_{q,c}^n, \mathbf{s}^{p,\pi}_c).
\end{align*}

By combining the optimal values of the scenario clusters, it results the  lower bound $z_{SMC}(e^*)=\sum_{c\in\C} \bar{W}_c \,  z_c(e^*)$ of optimal value $z^*$, where $\bar{W}_c=\sum_{\omega \in \Omega_c} \hat{w}^\omega$.
It is worth mentioning that for $e^*=E-1$ $SMC$ is equivalent to $SWS$.

\section{Computational experiments}\label{sec:results}

In this section we report computational results for a case study inspired by a building complex in southern Germany. Results are provided for three strategic scenario-tree sizes (small, medium, large) and for three model variants: the RES model without user discomfort, the risk-neutral RES model, and the risk-averse RES model.

The remainder of this section is organized as follows.
Section \ref{sec:data} describes the main features of the case study.
Section \ref{sec:res-results-cpx} reports the results obtained by solving the RES model as a monolithic problem.
Section \ref{sec:results-LB} presents the lower bounds derived from the schemes introduced in Section \ref{sec:lb}.
Section \ref{sec:res-results-heur} analyzes the performance of the proposed matheuristic and compares it with an established matheuristic from the literature.
Section \ref{sec:vsd} quantifies the value of accounting for uncertainty.
Finally, Section \ref{sec:Risk-analysis} examines the impact of risk-averse modeling.

\subsection{The data}
\label{sec:data}
The case study is constructed using publicly available data sources.
Electricity prices are obtained from the ENTSOE's Transparency Platform, see \cite{ENTSOE}, while PV production parameters and load profiles are derived from the CoSSMic European project, see \cite{OPSD}, which provides detailed consumption profiles, with resolutions reaching individual devices, for various small businesses and residential households.

The problem dimensions for these models are summarized in Table \ref{pro-res-dimensions}.
It gives the size of the strategic scenario tree, the structure of the operational subtrees, and the cardinalities of the main RES sets.
\textcolor{black}{In all instances, each strategic stage corresponds to one year of the planning horizon.
}
{\color{black}
Operational uncertainty within each stage is represented through a finite set of operational scenarios, each corresponding to a representative day. Representative days are selected by applying a k-medoids clustering algorithm to historical daily profiles of electricity demand and PV generation, following the methodology proposed in \cite{Micheli20}. Each cluster is represented by its medoid day, and its probability is set equal to the proportion of historical days assigned to that cluster. A detailed description of the clustering procedure is provided in Appendix C.
}

All three instances (small, medium and large) include investment in three PV technologies (say, mono-crystalline, poly-crystalline and thin-film panels) and two types of BESS (say, lead-acid and lithium-ion).
\begin{longtable}[ht!]
{l|rrrr|rr|rrrrrr}
\hline
Size & $|\E|$ & $|\N|$ & $|\Omega|$ & $\tilde{\S}_1$
 & $\tilde{\Pi}$ & $\tilde{\T}$
 & $|\I|$ & $|\BB|$ & $|\J_1|$ & $|\J_2|$ & $|\H_1|$ & $|\H_2|$
\\
\hline
Small  & 3 & 13 & 9 & 3 & 10 & 24 & 3 & 2 & 25 & 25 & 10 & 10 \\
Medium  & 4 & 40 & 27 & 3 & 20 & 24 & 3 & 2 & 40 & 35 & 15 & 15  \\
Large  & 6 & 364 & 243 & 3 & 20 & 24 & 3 & 2 & 75 & 75 & 50 & 50 \\
\hline
\caption{Dimensions of the strategic scenario tree ($|\E|$, $|\N|$, $|\Omega|$, $\tilde{\S}_1$), operational subtrees ($\tilde{\Pi}$, $\tilde{\T}$), and RES sets ($|\I|$, $|\BB|$, $|\J_1|$, $|\J_2|$, $|\H_1|$, $|\H_2|$) for the small, medium, and large instances. Note: $\tilde{\S}_1$ represents the branching degree of the strategic scenario tree.}
\label{pro-res-dimensions}
\end{longtable}
Mono-crystalline panels are the most efficient but also the most expensive, while thin-film panels offer lower investment costs at the expense of reduced electricity production capacity. Regarding BESS, lithium-ion batteries have higher costs but greater storage capacity and higher charge/discharge depth bounds.
At the start of the planning horizon, investment costs are as follows: PV mono-crystalline 2.5 \texteuro$/W$, PV poly-crystalline 2.1 \texteuro$/W$, PV thin-film 1.95 \texteuro$/W$, lead-acid BESS 1.05 \texteuro$/W$, and lithium-ion BESS 1.3 \texteuro$/W$.
In subsequent stages, investment costs in children nodes of each node $n$ are determined based on three trajectories: a stable cost trajectory, a cost reduction up to 30\%, and a cost increase up to a 30\% compared to the ancestor node $n$.
Maintenance costs are assumed to be 1.5\% of the investment costs in each strategic node.
A budget of 20,000 \texteuro$\,$ for RES infrastructure investments is imposed at each strategic node.
{\color{black}
At the operational level, the three instances share an hourly time discretization over a 24-hour horizon.
}

The instances differ in the number of controllable loads: 50 for the small instance, 75 for the medium instance, and 150 for the large instance.
To analyze the impact of discomfort modeling on the optimal solution of the problem, we evaluate three models for each instance type:
\begin{enumerate}
\item No-discomfort contention model \eqref{of}--\eqref{tac-battery_flow_balance}, from now on \textit{No-D}.
\item Risk-neutral model \eqref{of}--\eqref{rn}, from now on \textit{RN}, to bound the expected user discomfort.
\item Risk-averse model \eqref{of}--\eqref{sd}, from now on \textit{SD}, to prevent extreme operational scenarios.
\end{enumerate}
\textcolor{black}{
In model \textit{RN}, the maximum expected discomfort at each strategic node is bounded by the parameter $\hat{D}_e$, see constraint \eqref{rn}. In the computational experiments, this threshold is set to $\hat{D}_e = 20$ for the small- and medium-sized instances and to $\hat{D}_e = 40$ for the large-sized instance, reflecting the different number of households in the building complex.
In model \textit{SD}, a single policy profile is considered for each strategic node (i.e., $|\P_e| = 1$ for all $e \in \E$). The parameters of the stochastic dominance constraint system \eqref{sd} are fixed as follows in all experiments.
The probability bound is set to $\overline{\eta}^p = 0.05$, limiting the probability that the total operational discomfort exceeds the threshold to at most 5\%.
The second-order stochastic dominance constraint system impose a maximum excess fraction $\overline{s}^p = 0.25$ and a maximum expected excess fraction $\overline{\overline{s}}^p = 0.05$.
The discomfort threshold $\overline{D}^p$ coincides with the maximum expected discomfort $\hat{D}_e$ used in model \textit{RN}.
}

In total, nine models are tested.
Computations were performed on an ASUS laptop equipped with a 3 GHz Intel Core i7-5500U processor and 4 GB of RAM, using GAMS 48.2.0 and Gurobi 12.0.2.
We used GUSS (Gather-Update-Solver-Scatter) to process batches of structurally identical problems, reducing repeated model compilation overhead (see \cite{bussieck}).

\subsection{Results of the straightforward solution of the RES models}
\label{sec:res-results-cpx}

Table \ref{results-cpx} reports the dimensions, the solution cost and the CPU time for solving the nine models derived from applying the three versions of the monolithic RES MH model.
The headings are as follows: $\#Cons$, number of constraints; $\#Int$, number of general integer variables; $\#01$, number of binary variables; $\#Cont$, number of continuous variables; $\#Nze$, number of nonzero coefficients in the constraint matrix; and $z^*$ and $t^*$, optimal solution cost and CPU time (in seconds), respectively.
{\color{black}
\begin{longtable}
{ll|rrrrrrr}
\hline
Size & Model & $\#Cons$ & $\#Int$ & $\#01$ & $\#Cont$ & $\#Nze$ &
$z^*$ & $t^*$  \\
\hline
Small & \textit{No-D} & 317,386 & 26 & 93,650 & 93,756 & 625,435 & 4,375,801 & 38 \\
Small & \textit{RN} & 317,399 & 26 & 93,650 & 93,756 & 689,105 &
4,400,328 & 40 \\
Small & \textit{SD} & 317,685 & 26 & 93,780 & 93,886 & 744,025 & 5,526,312 & 51 \\ \hline
Medium & \textit{No-D} & 2,811,791 & 80 & 768,121 & 864,480 & 6,669,158 & 15,404,774 & 569
\\
Medium & \textit{RN} & 2,811,831 & 80 & 768,121 & 864,480 & 7,331,057 & 16,924,140 & 588
\\
Medium & \textit{SD} & 2,813,511 & 80 & 768,921 & 865,280 & 8,222,118 & 21,603,367 & 1151
\\ \hline
Large & \textit{No-D} & 43,181,342 & 728 & 12,232,240 & 10,490,300 & 96,961,232 & -- & OOM
\\
Large & \textit{RN} & 43,181,706 & 728 & 12,232,240 & 10,490,300 & 101,244,232
& -- & OOM
\\
Large & \textit{SD} &  43,197,036 & 728 & 12,239,520 & 10,497,580 &
121,689,557
& -- & OOM
\\
\hline
\caption{Models' dimensions, solution cost and CPU time. Note: OOM stands for Out Of Memory}
\label{results-cpx}
\end{longtable}
}
Small-sized instances are solved to optimality within seconds; medium-sized instances are solved within 20 minutes.
This performance is partly attributable to the strong formulation and problem structure.
In contrast, solving the large-sized models directly proves computationally infeasible, highlighting the need for alternative solution strategies.

For the instances that can be solved up to optimality, the results clearly show the economic impact of incorporating discomfort modeling.
In the small-sized instance, the risk-averse model \textit{SD} yields a total cost that is 26.30\% higher than the no-discomfort model \textit{No-D} and 25.59\% higher than the risk-neutral model \textit{RN}.
In the medium-sized instance, the cost increase is even more pronounced: the \textit{SD} solution is 40.24\% higher than \textit{No-D} and 27.64\% higher than \textit{RN}.
\textcolor{black}{
While these results quantify the additional cost associated with incorporating risk-averse measures in instances solvable up to optimality,
Section \ref{sec:Risk-analysis} analyzes
the corresponding benefits in terms of operational robustness and user discomfort mitigation for the large-sized instance.
}

\subsection{Results for the lower bound schemes}\label{sec:lb-results}
\label{sec:results-LB}
{\color{black}
This section presents lower bounds for the large-sized instance, for which direct solution of the RES model is computationally intractable.
The bounds are obtained using the schemes outlined in Section \ref{sec:lb} and are used to assess the quality of the solutions produced by the proposed matheuristic.
}

\begin{longtable} [ht!]
{l|rr|rr|rr}
\hline
Model & \multicolumn{2}{c|}{$SWS$} & \multicolumn{2}{c|}{$MHEV$} & \multicolumn{2}{c}{$MHOEV$}\\
\hline
\textit{No-D} & 20,497,013  & (1096)
& 2,206,149 & (81) & 2,348,221 & (290) \\
\textit{RN} & 23,110,918 & (1183)
& 2,897,976 & (86) & 2,901,853 & (293)  \\
\textit{SD} & 36,781,097 & (1279)
& \multicolumn{2}{c|}{$-$} & \multicolumn{2}{c}{$-$}  \\
\hline
\caption{Lower bound cost and CPU time (in parenthesis) in $SWS$, $MHEV$ and $MHOEV$ for the large instance}
\label{res-LB-b}
\end{longtable}

Table \ref{res-LB-b} shows $SWS$ yields much tighter bounds than $MHEV$ and $MHOEV$.
Those bounds can be further strengthened by grouping the strategic scenarios, as done in the $SMG$ scheme; the corresponding results are shown in Table \ref{res-LB-SMG} for different choices of the number of groups, namely $G=10, 20$.
Parameter $G$ introduces a trade-off between bound tightness and computational effort. Smaller values of $G$ yield tighter bounds but require solving larger subproblems, whereas larger values reduce computational burden at the expense of bound quality.
\begin{longtable}[ht!]
{l|rr|rr}
\hline
Model & \multicolumn{2}{c|}{$SMG(10)$} & \multicolumn{2}{c}{$SMG(20)$}\\
\hline
\textit{No-D} & 20,568,104 & (6642) & 20,647,634 & (7320) \\
\textit{RN} & 23,167,843 & (6782) & 23,235,098 & (7671) \\
\textit{SD} & 36,939,100 & (7009) & 37,002,098 & (8135) \\
 \hline
\caption{Lower bound cost and CPU time (in parenthesis) in $SMG$ for $G=10, 20$ for the large instance}
\label{res-LB-SMG}\\
\end{longtable}
The $SMC$ scheme further strengthens the $SWS$ lower bound by relaxing non-anticipativity constraints up to a breaking stage $e^*$.
Table \ref{res-LB-scm} shows the tight lower bounds obtained using $SMC$ for three distinct values of the breaking stage $e^*$.
The case $e^*=5$ is excluded from Table \ref{res-LB-scm}, since it corresponds to $SWS$, while for $e^*=1$, the corresponding model is not computationally tractable.
The resulting bounds improve as the breaking stage $e^*$ decreases, although the marginal gains diminish for smaller values of $e^*$.
This suggests that significant benefits are already captured at higher values of $e^*$.
From a computational perspective, decreasing $e^*$ reduces the number of submodels but increases their individual size, resulting in a non-monotonic behavior of CPU times.
The lowest CPU time is achieved for $e^* = 3$.

\begin{longtable}[ht!]
{l|rr|rr|rr}
\hline
Model & \multicolumn{2}{c|}{$SMC(4)$} & \multicolumn{2}{c|}{$SMC(3)$} & \multicolumn{2}{c}{$SMC(2)$} \\
\hline
\textit{No-D} & 20,656,110 & (3527) & 20,897,013 & (3229) & 20,902,001 & (6480)
 \\
\textit{RN} & 23,426,109 & (4116) & 23,621,018 & (3851) & 23,720,923 & (6721)
 \\
\textit{SD} & 37,022,913 & (5167) & 37,078,181 & (4204) & 37,150,227 & (7650)
 \\
 \hline
\caption{Lower bound cost and CPU time (in parenthesis) in $SMC$ with breaking stage $e^*= 2, 3, 4 $ for the large instance}
\label{res-LB-scm}\\
\end{longtable}

Several observations emerge from the lower-bound analysis.
First, the weak bounds obtained through expectation-based schemes confirm that ignoring short-term variability leads to significant cost underestimation.
Incorporating operational uncertainty is therefore essential to obtain meaningful lower bounds for the RES problem.
Second, $SWS$ provides substantially tighter bounds, which can be further strengthened through $SMG$ and $SMC$.
Among these, $SMC$ consistently yields the tightest bounds while maintaining moderate computational effort, thereby representing the most effective strengthening scheme for the large-sized instance in any of the three model variants.

\subsection{Performance assessment of SFR3}\label{sec:res-results-heur}

This section evaluates the performance of the proposed matheuristic SFR3 on models \textit{No-D}, \textit{RN}, and \textit{SD}. We first assess solution quality and then benchmark SFR3 against an alternative horizon-decomposition matheuristic.

\subsubsection{Solution quality of SFR3}\label{sec:res-results-sfr3}
Tables \ref{res-sfr3} and \ref{res-sfr3-a} report the parameter configuration $(\hat{e}, \hat{e}^R, \phi_e)$, the reference objective function value ($z^*$ for small and medium instances, and $\underline{z}$ for the large instance), the feasible solution cost $z_{SFR3}$, the corresponding optimality gap, and the CPU time $t_{SFR3}$.
\begin{longtable}[ht!]
{ll|rrr|rr|r|r}
\hline
Size & Model & $\hat{e}$ & $\hat{e}^R$ & $\phi_e$ & ${z^*}$ & $z_{SFR3}$ & $GAP(z^*)$ & $t_{SFR3}$  \\
\hline
Small & \textit{No-D} & 1 & 0 & 0 & 4,375,801 & 4,486,224 & 2.52\% & 36
\\

Small & \textit{No-D} & 2 & 0 & 0 & 4,375,801 & 4,377,674 & 0.04\% & 40
\\

Small & \textit{No-D} & 1 & 2 & $\frac{1}{3}$ & 4,375,801 & 4,376,737 & 0.02\% & 39
\\

Small & \textit{No-D} & 2 & 1 & $\frac{1}{3}$ & 4,375,801 & 4,375,907 & 0.00\% & 42
\\
\hdashline

Small & \textit{RN} & 1 & 0 & 0 & 4,400,328 & 4,512,602 & 2.55\% & 40
\\

Small & \textit{RN} & 2 & 0 & 0 & 4,400,328 & 4,488,247 & 2.00\% & 44
\\

Small & \textit{RN} & 1 & 2 & $\frac{1}{3}$ & 4,400,328 & 4,407,861
 & 0.17\% & 47
\\

Small & \textit{RN} & 2 & 1 & $\frac{1}{3}$ & 4,400,328 & 4,401,605
 & 0.03\% & 51
\\
\hdashline

Small & \textit{SD} & 1 & 0 & 0 & 5,526,312 & 5,646,741 & 2.18\% & 45
\\

Small & \textit{SD} & 2 & 0 & 0 & 5,526,312 & 5,549,202 & 0.41\% & 50
\\

Small & \textit{SD} & 1 & 2 & $\frac{1}{3}$ & 5,526,312 & 5,529,561 & 0.06\% & 52
\\

Small & \textit{SD} & 2 & 1 & $\frac{1}{3}$ & 5,526,312 & 5,527,876 & 0.03\% & 62
\\
\hline

Medium & \textit{No-D} & 1 & 0 & 0 & 15,404,774 & 15,989,385 & 3.79\% & 126
\\

Medium & \textit{No-D} & 2 & 0 & 0 & 15,404,774 & 15,654,316 & 1.62\% & 211
\\

Medium & \textit{No-D} & 1 & 2 & $\frac{1}{3}$ & 15,404,774 & 15,414,479
 & 0.06\% & 205
\\

Medium & \textit{No-D} & 2 & 1 & $\frac{1}{3}$ & 15,404,774 & 15,413,310
 & 0.06\% & 234
\\
\hdashline

Medium & \textit{RN} & 1 & 0 & 0 & 16,924,140 & 17,627,372 & 4.16\% & 214
\\

Medium & \textit{RN} & 2 & 0 & 0 & 16,924,140 & 17,434,234 & 3.01\% & 260
\\

Medium & \textit{RN} & 1 & 2 & $\frac{1}{3}$ & 16,924,140 & 17,028,703 & 0.61\% & 281
\\

Medium & \textit{RN} & 2 & 1 & $\frac{1}{3}$ & 16,924,140 & 17,000,389 & 0.45\% & 385
\\
\hdashline

Medium & \textit{SD} & 1 & 0 & 0 & 21,603,367 & 22,382,320 & 3.61\% & 269
\\

Medium & \textit{SD} & 2 & 0 & 0 & 21,603,367 & 22,050,665 & 2.07\% & 270
\\

Medium & \textit{SD} & 1 & 2 & $\frac{1}{3}$ & 21,603,367 & 21,800,606 & 0.91\% & 380
\\

Medium & \textit{SD} & 2 & 1 & $\frac{1}{3}$ & 21,603,367 & 21,633,381 & 0.14\% & 323
\\
\hline
\caption{Results of matheuristic SFR3 for obtaining feasible solutions for small- and medium-sized instances}
\label{res-sfr3}
\end{longtable}
Table \ref{res-sfr3} shows that the weak myopic strategy ($\hat{e}=1$, $\hat{e}^R=0$) yields poor solution quality across all models.
Increasing the look-ahead horizon to $\hat{e}=2$ improves performance for small instances but remains insufficient for medium-sized problems.
Strategies incorporating relaxation stages and probabilistic node selection achieve near-optimal solutions while maintaining moderate computational effort.
In particular, the strategy ($\hat{e}=2$,  $\hat{e}^R=1$) consistently outperforms the strategy ($\hat{e}=1$, $\hat{e}^R=2$), offering superior solution quality while requiring comparable computational effort, indicating that balanced horizon depth and controlled relaxation improve both accuracy and efficiency.

In the large-sized instance, it can be observed in Table \ref{res-sfr3-a} that the weak and stronger myopic strategies produce weak solutions.
Increasing the horizon to $\hat{e}=3$ yields moderate improvements but at a substantial computational cost.
Accurate solutions require incorporating relaxation stages with $\phi_e = \frac{1}{3}$.
The strategies ($\hat{e}=2$, $\hat{e}^R=2$) and ($\hat{e}=3$, $\hat{e}^R=1$) achieve optimality gaps below 0.32\% for all the three models.
Observe that the strategy ($\hat{e}=2$, $\hat{e}^R=2$) provides the best balance, delivering comparable solution quality with lower computational effort.

\begin{longtable}
{ll|rrr|rr|r|r}
\hline
Size & Model & $\hat{e}$ & $\hat{e}^R$ & $\phi_e$ & $\underline{z}$ & $z_{SFR3}$ & $GAP(\underline{z})$ & $t_{SFR3}$ \\
\hline
Large & \textit{No-D} & 3 & 1 & 0 & 20,902,001 & 21,954,417 & 5.03\% & 1470
\\
Large & \textit{No-D} & 2 & 0 & 0 & 20,902,001 & 21,468,821  & 2.71\%  &  51785
\\
Large & \textit{No-D} & 3 & 0 & 0 & 20,902,001 & 21,218,478 & 1.51\% & 7020
\\
Large & \textit{No-D} & 2 & 2 & $\frac{1}{3}$ & 20,902,001 & 20,919,391 & 0.08\%  &  2367
\\
Large & \textit{No-D} & 3 & 1 & $\frac{1}{3}$ & 20,902,001 & 20,912,766 & 0.05\% &  8963
\\
\hdashline
Large & \textit{RN} & 1 & 0 & 0 & 23,720,876 &  24,672,083 & 4.01\% &  1574
\\
Large & \textit{RN} & 2 & 0 & 0 & 23,720,876 &  24,483,170 & 3.21\% &   1869
\\
Large & \textit{RN} & 3 & 0 & 0 & 23,720,876 &  24,263,610 & 2.29\% &   6895
\\
Large & \textit{RN} & 2 & 2 & $\frac{1}{3}$ & 23,720,876 &  23,744,882 & 0.10\%  &  2467
\\
Large & \textit{RN} & 3 & 1 & $\frac{1}{3}$ & 23,720,876 &  23,747,230 & 0.11\% &  9063
\\
\hdashline
Large & \textit{SD} & 1 & 0 & 0 & 37,150,197 & 39,310,109 & 5.81\% &  1710
\\
Large & \textit{SD} & 2 & 0 & 0 & 37,150,197 & 38,520,593 & 3.69\% &  1869
\\
Large & \textit{SD} & 3 & 0 & 0 & 37,150,197 & 38,100,425 & 2.56\% &  7778
\\
Large & \textit{SD} & 2 & 2 & $\frac{1}{3}$ & 37,150,197 & 37,268,929 & 0.32\% &  2369
\\
Large & \textit{SD} & 3 & 1 & $\frac{1}{3}$ & 37,150,197 & 37,268,520 & 0.32\% & 9294
\\
\hline
\caption{Results of matheuristic algorithm SFR3 for obtaining feasible solutions for large-sized models}
\label{res-sfr3-a}
\end{longtable}


\bcbl
\subsubsection{Benchmark comparison with matheuristic SRH}\label{sec:res-results-srh}
To further assess the performance of SFR3, we benchmark it against matheuriatic SRH (it stands for shrinking-rolling horizon), a well-established horizon-decomposition approach in stochastic optimization introduced in \cite{BG04}. SRH approximates a multistage stochastic program by solving a sequence of smaller two-stage stochastic problems. At each iteration, decisions corresponding to the earliest stage are fixed and implemented, while the remaining uncertainty is progressively resolved as the horizon shrinks.
As a result, SRH partially mitigates the myopic behavior typical of rolling-horizon methods, at the cost of solving a larger number of stochastic subproblems.

Table~\ref{res-srh} reports the results obtained by solving the test instances with SRH and compares them with those produced by SFR3. For small- and medium-sized instances, $\hat{z}$ denotes the optimal solution cost $z^*$, while for the large-sized instance it denotes the tightest available lower bound $\underline{z}$. On the other hand, $z_{\mathrm{SRH}}$ represents the cost of the feasible solution obtained by SRH, and the corresponding optimality gap is computed as $(z_{\mathrm{SRH}} - \hat{z}) / \hat{z}$. The computational effort required by SRH is reported as $t_{\mathrm{SRH}}$.
To benchmark SFR3 against SRH, the table also reports the goodness ratio $\mathrm{GR}_{\mathrm{SFR3}} = z_{\mathrm{SFR3}} / z_{\mathrm{SRH}}$ and the CPU time ratio $\mathrm{TR}_{\mathrm{SFR3}} = t_{\mathrm{SFR3}} / t_{\mathrm{SRH}}$. Values of $\mathrm{GR}_{\mathrm{SFR3}} < 1$ indicate that SFR3 achieves a lower solution cost than SRH, while values of $\mathrm{TR}_{\mathrm{SFR3}} < 1$ indicate a reduction in computational time.

\begin{longtable}
{ll|rr|r|rrr}
\hline
Size & Model  & $\hat{z}$ & $z_{\mathrm{SRH}}$ & $GAP$ & $t_{\mathrm{SRH}}$ & $GR_{SFR3}$ & $TR_{SFR3}$\\
\hline
Small & \textit{No-D} & 4,375,801 & 4,380,548 & 0.11\% & 44 & 0.999 & 0.955 \\

Small & \textit{RN} & 4,400,328 & 4,435,461 & 0.80\% & 55 & 0.992 & 0.927
\\
Small & \textit{SD} & 5,526,312 & 5,534,688 & 0.15\% & 70 & 0.999 & 0.827
\\
\hdashline
Medium & \textit{No-D} & 15,404,774 & 15,489,881 & 0.55\% & 521 & 0.995 & 0.449 \\

Medium & \textit{RN} & 16,924,140 & 16,987,240 & 0.32\% & 661 & 1.001 & 0.582
\\
Medium & \textit{SD} & 21,603,367 & 21,646,212 & 0.19\% & 821 & 0.999 & 0.393
\\
\hdashline
Large & \textit{No-D} & 20,902,001 & 21,285,417 & 1.83\% & 20,301 & 0.982 & 0.442 \\

Large & \textit{RN} & 23,720,876 & 23,991,631 & 1.14\% & 21,496 & 0.990 & 0.115
\\
Large & \textit{SD} & 37,150,197 & 37,545,644 & 1.06\% & 22,082 & 0.993 & 0.421
\\
\hline
\caption{Results of matheuristic SRH for obtaining feasible solutions for the instances}
\label{res-srh}
\end{longtable}

{\color{black}
For the small-sized instances, SFR3 and SRH achieve comparable solution quality, with negligible differences in optimality gaps and similar computational times.
For the medium-sized instances, SFR3 generally attains a slightly smaller solution cost while requiring substantially less CPU time.
The only exception is the \textit{RN} case, where SRH provides a marginally smaller cost;
however, this improvement is obtained at a significantly higher computational cost.
For the large-sized instances, SFR3 consistently delivers smaller solution cost than SRH and reduces CPU times dramatically,
see the related $GR_{SFR3}$ and $TR_{SFR3}$, confirming that SFR3 scales more effectively as problem size increases.
Overall, these results demonstrate that SFR3 provides a superior balance between solution accuracy and computational efficiency compared to the shrinking-horizon benchmark.
}

\subsection{Value of stochastic modeling}
\label{sec:vsd}
Quantifying the benefit of stochastic modeling in multistage settings has received limited attention, and existing evidence is largely application-specific (e.g., supply network management in \cite{vss07}).
In the multi-horizon context, \cite{Maggioni20} provides, to the best of our knowledge, the first systematic discussion and quantification framework.
Accordingly, to assess the benefit of explicitly accounting for uncertainty in strategic decisions within our multistage multi-horizon model, we report the Value of Strategic Decision (VSD) following \cite{Maggioni20}.
The VSD measures the loss incurred when strategic decisions are obtained from an expected-value model and then evaluated under uncertainty.
In our setting, it is computed by comparing two feasible solutions of the full multi-horizon stochastic problem.
As stochastic benchmark, we use the high-quality feasible solution returned by SFR3, which is near-optimal for the large-sized instance.
The comparison solution is obtained by fixing the strategic decision variables to the values retrieved from the solution of the expected-value model $MHEV$ and re-optimizing all operational decisions in model \textit{SD}.
This procedure preserves full operational recourse and yields a feasible upper bound, denoted by $z_{S-MHEV}$, on the stochastic cost associated with the expectation-based strategic decisions.
The VSD is defined as the difference between cost $z_{S-MHEV}$ and cost $z_{SFR3}$ (taken from Table \ref{res-sfr3-a}), so, VSD should be interpreted as a conservative estimate, i.e., a lower bound on the true value of stochastic modeling.
Table~\ref{tab:vsd} reports the resulting VSD for the large-sized instance;
it is strictly positive for all three model variants, indicating that the stochastic strategic decisions derived from model $MHEV$ lead to systematically higher costs when evaluated under uncertainty, even after allowing for full operational recourse.
Measures that fix both strategic and operational decisions (e.g., the Value of the Stochastic Operational Decision, VSOD) are not reported, as they yield infeasible solutions in this application, consistently with the observations in \cite{Maggioni20}.
\begin{longtable}[ht!]
{lrrrrr}
\hline
Model
& $z_{S-MHEV}$ [\texteuro]
& CPU time [s]
&
$z_{SFR3}$
& VSD [\texteuro]
& $GR$ \\
\hline
\textit{No-D}  & 28,260,495 & 13,188 & 20,912,766 &  7,347,729 & 0.74 \\
\textit{RN}    & 34,922,397 & 15,022 & 23,747,230 & 11,175,167 & 0.68 \\
\textit{SD}    & 59,156,790 & 17,901 & 37,268,929 & 21,887,861 & 0.63 \\
\hline
\caption{Value of Strategic Decision (VSD) for the large-sized instance
%
%
}
\label{tab:vsd}
\end{longtable}
The VSD magnitude increases markedly when moving from model \textit{No-D} to models \textit{RN} and \textit{SD}.
In particular, VSD exceeds 21 M\texteuro $\,$ in model \textit{SD}, highlighting that ignoring uncertainty at the strategic level becomes increasingly costly as discomfort and risk considerations are strengthened.
To further contextualize these results, the goodness ratio $GR$, to be expressed as $z_{SFR3}/z_{S-MHEV}$, indicates an improvement achieved by explicitly accounting for uncertainty at the strategic level, for $GR<1$.
The decreasing goodness ratio observed across the three model variants confirms that the relative benefit of stochastic modeling becomes more pronounced as the problem places greater emphasis on discomfort and risk mitigation.
Overall, these results provide quantitative validation of the proposed multi-horizon stochastic framework.

\subsection{Impact of risk-averse modeling}
\label{sec:Risk-analysis}

We focus on the large-size instance, which represents the most challenging setting in terms of uncertainty and computational complexity, to assess how risk-averse modeling affects strategic investment decisions and short-term operational behavior.
The analysis is based on the SFR3 solution and compares models \textit{No-D}, \textit{RN}, and \textit{SD} under the parameter settings given in Section~\ref{sec:data}. Table~\ref{tab:capacity-stages} reports the evolution of installed PV and BESS capacities across the stages.
\begin{longtable}[ht!]
{llrrrrrr}
\hline
Model & Technology & e=1 & e=2 & e=3 & e=4 & e=5 & e=6 \\
\hline
\multirow{3}{*}{\textit{No-D}}
& Polycrystalline PV
& 10.00 & 15.00 & 18.66 & 20.27 & 20.55 & 20.55 \\
&Lead-acid BESS
& 0.00 & 0.00 & 0.55 & 1.55 & 1.66 & 1.66 \\
&Lithium-ion BESS
& 4.00 & 10.00 & 13.66 & 15.85 & 18.20 & 18.20 \\
\hline
\multirow{3}{*}{\textit{RN}}
&Polycrystalline PV
& 6.00 & 11.33 & 17.39 & 22.27 & 23.05 & 23.68 \\
&Lead-acid BESS
& 0.00 & 1.33 & 3.00 & 4.81 & 4.81 & 4.81 \\
&Lithium-ion BESS
& 9.00 & 16.00 & 19.65 & 19.85 & 20.00 & 20.00 \\
\hline
\multirow{3}{*}{\textit{SD}}
&Polycrystalline PV
& 4.00 & 9.30 & 14.41 & 20.82 & 23.21 & 23.78 \\
&Lead-acid BESS
& 0.00 & 5.00 & 5.00 & 5.00 & 5.00 & 5.00 \\
&Lithium-ion BESS
& 16.00 & 18.00 & 19.78 & 20.00 & 20.00 & 20.00 \\
\hline
\caption{Installed power capacity (kW) \textit{by} stage for the large-sized instance under the three model variants.
Reported PV capacities refer to polycrystalline panels; other PV technologies are not selected in the optimal solutions
}
\label{tab:capacity-stages}
\end{longtable}

Clear differences emerge across the three models. Investment decisions in model \textit{No-D} primarily favor cost-reducing generation capacity, with limited emphasis on flexibility-oriented assets.
Introducing the risk-neutral bound \eqref{rn} leads to a partial investment re-balancing by increased deployment of battery energy storage to mitigate average operational stress.
This shift becomes more pronounced in model \textit{SD}, where investment decisions place greater emphasis on battery energy storage in the early stages of the planning horizon.
At the same time, PV capacity expansion in \textit{SD} is more gradual, with lower installation levels in the initial stages and more pronounced increases toward the terminal one.
This behavior is consistent with the presence of binding budget constraints, which limit the possibility of simultaneously expanding all technologies at the outset and induce a prioritization of assets with higher short-term risk-mitigation value.
Table \ref{tab:energy-balance} summarizes the expected contribution of the RES main components to the daily energy balance at the planning horizon in the three model variants.
\begin{longtable}[ht!]
{lrrrrrr}
\hline
Model & PV & Import & BESS Discharge & Load curtailment & BESS Charge & Total Load \\
\hline
\textit{No-D} & 73.29  & 141.32 & 45.21 & 54.97 & 46.42 & 268.37 \\
\textit{RN}   & 81.30  & 172.79 & 50.78 & 15.54 & 52.04 & 268.37 \\
\textit{SD}   & 82.01  & 182.21 & 55.21 &  5.35 & 56.41 & 268.37 \\
\hline
\caption{End-of-horizon (i.e., stage 6) expected daily energy balance (kWh/day) for the large-sized instance}
\label{tab:energy-balance}
\end{longtable}
Notice that a substantial load curtailment is performed in model \textit{No-D}, while storage is used primarily for limited arbitrage.
Introducing the risk-neutral bound \eqref{rn} on expected discomfort significantly alters this behavior;
observe that curtailment is markedly reduced, from about 55 kWh/day to approximately 16 kWh/day, and the system increasingly relies on a combination of grid imports, storage operation, and PV generation to supply demand.
However, despite the higher use of flexibility in the resources on average, non-negligible curtailment remains, indicating that expectation-based control alone does not fully prevent stressed operational outcomes.
Additionally, introducing the risk-averse constraint system \eqref{sd} further reshapes the operational balance;
observe that curtailment is reduced to a marginal level, while battery charging and discharging increase substantially, reflecting a systematic reliance on storage-based flexibility to hedge against adverse realizations.
Compared to the \textit{RN} solution, the \textit{SD} one exhibits higher storage throughput and slightly higher grid imports,
but achieves a much smoother supply structure which drastically reduces unmet demand.

\begin{figure}[ht!]
    \centering
    \includegraphics[width=0.65\linewidth]{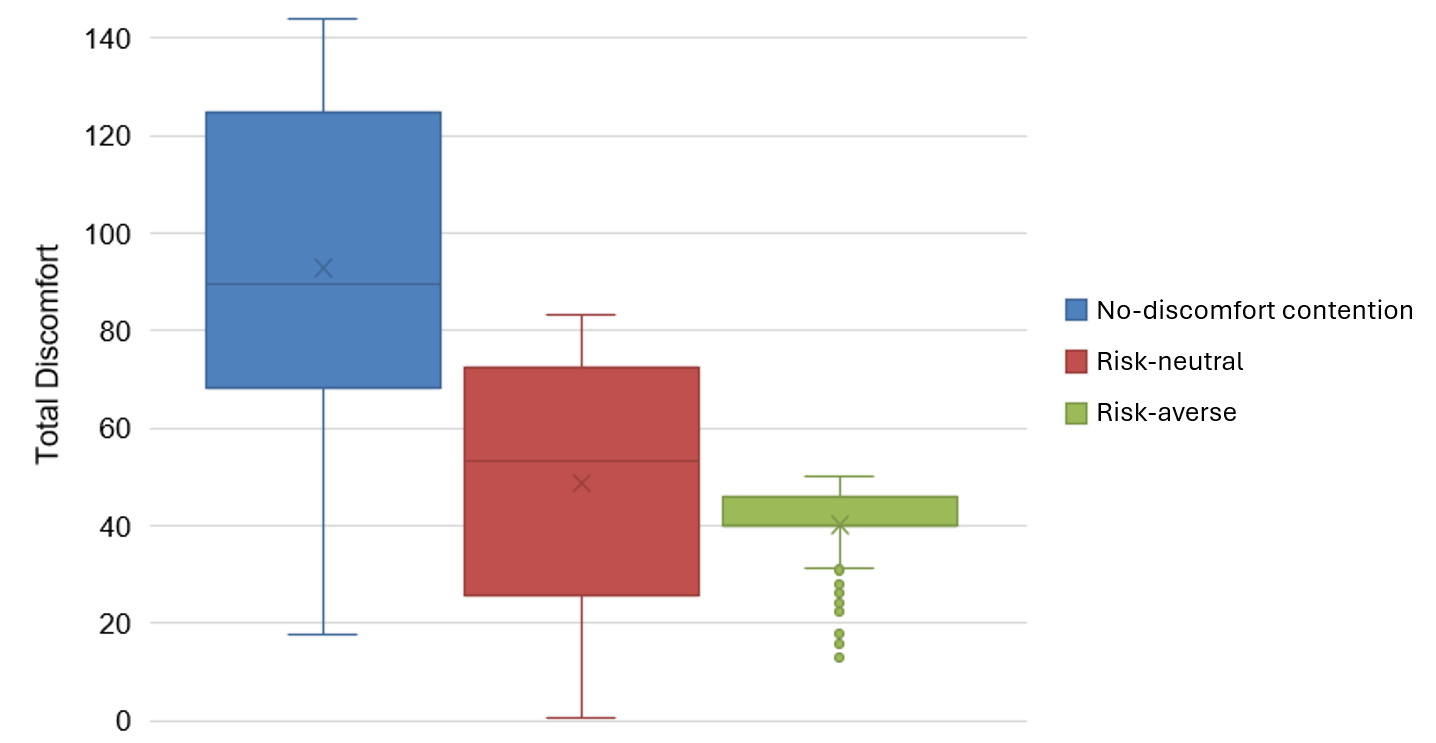}
    \caption{Distribution of nodal discomfort for the large-sized instance under the no-discomfort contention, the risk-neutral, and the risk-averse formulations.
    Note: The central boxes span from the 1st quantile to the 3rd one}
     \label{fig:discomfort-boxplot}
\end{figure}

Figure \ref{fig:discomfort-boxplot} illustrates the distribution of nodal discomfort value and its accumulated frequency across the strategic nodes and operational scenarios for the three model variants.
Model \textit{No-D} in the figure exhibits a wide dispersion of nodal outcomes,
with several operational realizations characterized by very large discomfort values.
Introducing the risk-neutral bound \eqref{rn} on expected discomfort substantially reduces average nodal discomfort;
however, severe nodal discomfort can still occur in a non-negligible number of operational scenarios.
Additionally, introducing the risk-averse constraint system \eqref{sd} leads to a qualitatively different behavior.
As shown in the figure, the upper tail of the nodal discomfort distribution becomes markedly compressed,
indicating that extreme nodal outcomes are effectively mitigated.

The risk-averse effect shown in Fig. \ref{fig:discomfort-boxplot} is further quantified in Table \ref{tab:discomfort-summary},
where it can be observed that the discomfort exceedance in the strategic nodes remains admissible, their occurrence is infrequent and their expected magnitude is substantially reduced, in line with the constraint system \eqref{sd} imposed in the model.
Additionally, observe that the risk-averse solution exhibits a slightly smaller mean nodal expected discomfort than the risk-neutral one,
indicating that the additional cost, see Table \ref{tab:vsd}, is due to structural changes in the allocation of resources aimed at mitigating extreme outcomes.
The benefit of the risk aversion becomes particularly evident when examining tail-related indicators.
Notice that the 95th percentile of nodal expected discomfort is reduced from 79.2 in model \textit{RN} and to 50.0 in model \textit{SD},
showing that extreme but plausible discomfort realizations are effectively capped.
This effect is reinforced by the sharp decrease in both the nodal threshold's maximum exceedance and violation frequency.
See that the latter decreases from more than 40\% of operational scenarios on average in model \textit{RN} to below 5\% in model \textit{SD}.
On the other hand, the maximum violation frequency across nodes is reduced to the 5\% prescribed bound,
indicating that even the most exposed nodes are effectively protected against recurrent high-discomfort realizations.

\begin{longtable}[ht!]
{lrrr}
\hline
& \textit{No-D} &  \textit{RN}  &  \textit{SD} \\
\hline
Mean nodal expected discomfort in LHS \eqref{sd-1}                       & 86.93 & 39.94 & 36.65 \\
95th percentile nodal expected discomfort in LHS \eqref{sd-1}            & 139.7 & 79.2 & 50.0 \\
Mean nodal maximum discomfort threshold's exceedance in LHS \eqref{sd-2}  & 58.22 & 27.08 & 9.35 \\
Mean nodal violation frequency in \eqref{sd_3bis}                        & 69.27\% & 40.54\% & 4.79\% \\
Maximum nodal violation frequency in \eqref{sd_3bis}                     & 84.67\% & 60.65\% & 5.00\% \\
\hline
\caption{Statistics of nodal discomfort for the large-sized instance in the three model variants. Note: violation occurs when the discomfort value is higher than the threshold}
\label{tab:discomfort-summary}
\end{longtable}

\bcb

\section{Conclusions and future research agenda}
\label{sec:conclu}


Motivated by the urgent need to accelerate the transition to sustainable energy systems, in this we address the challenge of designing a domestic renewable energy system, by proposing a MH optimization model that co-optimizes investment and operations under short- and long-term uncertainty,
and that embeds first- and second-order stochastic dominance constraints to mitigate extreme operational discomfort.
The proposed model is particularly well-suited to domestic RES design for three key reasons.
First, it includes long-term uncertainty in investment costs, enabling the designer to identify the optimal mix of PV and BESS technologies for the RES infrastructure over an extended horizon.
Second, the model considers a detailed representation of short-term operations, crucial for accurately simulating the actual functioning of the system.
Only by introducing in the problem an accurate short-term modeling, the challenges associated with the management of the intermittent solar generation and the supply of different types of loads can be properly addressed.
Third, the risk-averse structure limits exposure to extreme operational scenarios, enhancing system reliability and user protection.
Due to the significant computational complexity of MH stochastic models, the optimal solution for large-sized instances may be computationally intractable using state-of-the-art solvers.

\newpage
To overcome this challenge, we propose the matheuristic SFR3, specifically tailored to the structural properties of the proposed formulation.
\textcolor{black}{The algorithm generates high-quality feasible solutions efficiently and is benchmarked against an established shrinking-horizon matheuristic, consistently outperforming it in terms of scalability and overall performance.}
In addition, several lower-bounding schemes are implemented to rigorously assess solution quality.

The proposed methodology is validated using a case study based on a building complex in South Germany.
We investigate a range of instance sizes and model variants to rigorously evaluate the performance of the proposed method across different system configurations.
\bcbl
Numerical results demonstrate that SFR3 is capable of generating near-optimal solutions, with optimality gaps as low as 0.32\%, within a reasonable computational time, even for instances that are computationally intractable in monolithic form.
Furthermore, the Value of Strategic Decision analysis confirms the significant benefit of modeling uncertainty, while the risk analysis quantifies the protection provided by SD-based discomfort control.
These results indicate that the proposed framework provides a computationally tractable and practically relevant approach for domestic RES planning under multiscale uncertainty, particularly in settings where direct monolithic optimization is infeasible.

As future research, we aim to investigate the integration of advanced decomposition techniques, such as Stochastic Dual Dynamic Integer Programming (SDDiP), see \cite{Zou19}.
In particular, the matheuristic SFR3 could be employed to generate high-quality initial feasible solutions within an SDDiP framework that is adapted to tactical multi-horizon problems, potentially enhancing convergence behavior and enabling further scalability.
\bcb

\vspace{.2cm}
\bcbl
{\bf CRedIT author statement:}
\textbf{G. Micheli}: Conceptualization, Methodology, Formal analysis, Investigation, Data curator, Software, Writing-reviewing and editing;
\textbf{L.F. Escudero}: Conceptualization, Methodology, Investigation, Writing original draft preparation, Writing-reviewing and editing, Supervision;
\textbf{F. Maggioni}: Methodology, Investigation, Writing-reviewing and editing, Validation;
\textbf{G. Bayraksan}: Writing-reviewing and editing.
\bcb

\vspace{.2cm}
{\bf Acknowledgements:}
G. Micheli acknowledges the support from Gruppo Nazionale per il Calcolo Scientifico (GNCS-INdAM).
L.F. Escudero acknowledges the support from the projects RTI2018-094269-B-I00 and PID2021-122640OB-I00. \
F. Maggioni acknowledges the support from the PRIN2020 project "ULTRA OPTYMAL - Urban Logistics
and sustainable TRAnsportation: OPtimization under uncertainTY and MAchine Learning", funded by the Italian
University and Research Ministry (Grant 20207C8T9M,  \url{https://ultraoptymal.unibg.it}) and from Gruppo Nazionale per il Calcolo Scientifico (GNCS-INdAM).
G. Bayraksan acknowledges support from the U.S.\ Department of Energy, Office of Science, Office of Advanced Scientific Computing Research  under Grant DE-SC0023361.

\vspace{.2cm}
{\bf Declarations}: The authors have no competing interests to declare that are relevant to the content of this article. The data that support the findings of this study are openly available at \url{https://data.open-power-system-data.org/household_data} and \url{https://transparency.entsoe.eu/dashboard/show}, and upon request to the authors.

\end{document}